\documentclass[preprint,12pt]{elsarticle}

\usepackage[top=1.5cm, bottom=2.0cm, outer=2.25cm, inner=2.25cm, heightrounded, marginparwidth=2.5cm, marginparsep=0.5cm]{geometry}

\usepackage[english]{babel}
\usepackage[latin1]{inputenc}
\usepackage[T1]{fontenc}
\usepackage{graphicx}
\usepackage{booktabs}
\usepackage{caption}
\usepackage[labelformat=simple]{subcaption}

\usepackage{color}

\usepackage{amsmath}
\usepackage{amsfonts}
\usepackage{amssymb}
\usepackage{amsthm}
\usepackage[ruled]{algorithm}
\usepackage{algorithmic}
\newcommand{\R}{\mathbb{R}}  
\newcommand{\N}{\mathbb{N}}  
\newcommand{\E}{\mathbb{E}}  
\newcommand{\norm}[1]{\left\Vert#1\right\Vert}
\DeclareMathOperator*{\argmin}{arg\,min}
\DeclareMathOperator{\bias}{bias}
\DeclareMathOperator{\var}{Var}
\DeclareMathOperator{\cov}{Cov}
\begin{document}

\begin{frontmatter}
  \title{A new framework for extracting coarse-grained models from time series with
    multiscale structure}
%
\author[CE]{S.~Kalliadasis}
\ead{s.kalliadasis@imperial.ac.uk}
\author[MA]{S.~Krumscheid\corref{cor1}}
\ead{s.krumscheid10@imperial.ac.uk}
\author[MA]{G.~A.~Pavliotis}
\ead{g.pavliotis@imperial.ac.uk}
\cortext[cor1]{Corresponding author}
%
\address[CE]{Department of Chemical Engineering, Imperial
  College London, London SW7 2AZ, UK}
\address[MA]{Department of Mathematics, Imperial
  College London, London SW7 2AZ, UK}
%
\begin{abstract}
  In many applications it is desirable to infer coarse-grained models
  from observational data. The observed process often corresponds only
  to a few selected degrees of freedom of a high-dimensional dynamical
  system with multiple time scales. In this work we consider the
  inference problem of identifying an appropriate coarse-grained model
  from a single time series of a multiscale system. It is known that
  estimators such as the maximum likelihood estimator or the quadratic
  variation of the path estimator can be strongly biased in this
  setting. Here we present a novel parametric inference methodology
  for problems with linear parameter dependency that does not suffer
  from this drawback. Furthermore, we demonstrate through a wide
  spectrum of examples that our methodology can be used to derive
  appropriate coarse-grained models from time series of partial
  observations of a multiscale system in an effective and systematic
  fashion.
\end{abstract}


\end{frontmatter}


%
%
\section{Introduction}
\label{sec:intro}

Many natural phenomena and technological applications are
characterized by the presence of processes occurring across different
length and/or time scales. Examples range from biological systems
\cite{Chauviere2010} and problems in atmosphere and ocean sciences
\cite{Majda2008,Culina2010} to molecular dynamics \cite{Griebel2007},
materials science \cite{Fish2009} and fluid and solid mechanics
\cite{Huerre1998,Horstemeyer2010,Savva2010}, to name but a
few. Studying the full dynamics of such systems is often a very
intricate task due to the complex structure of the systems which also
hampers the ability to obtain governing equations from first
principles. However, it is often possible to exploit, e.g., scale
separation in order to obtain a reduced (low-dimensional) model for a
few selected degrees of freedom. The coefficients and/or parameters in
the reduced model must be derived from the full dynamics through an
appropriate coarse-graining procedure; see, e.g.\ \cite{Chorin2000,
  Turkington2013, Venturi2014} for recent works on various
coarse-graining methods. As mentioned above, it is often not possible
to obtain such a coarse-grained equation in explicit form and one must
necessarily resort to observations
\cite{Pavliotis2007,Krumscheid2013}. It is thus desirable to
appropriately fit a reduced stochastic coarse-grained model to the
observations of the underlying complex process.

The general problem of obtaining a reduced coarse-grained model from
the full system can be formulated as follows. Let the underlying
system be given in terms of a dynamical system $Z$ which evolves,
symbolically written, according to the dynamics
\begin{equation}
  \frac{dZ}{dt} = F(Z)\;,\label{eq:dyn:sys:prototype}
\end{equation}
where the state space $\mathcal{Z}$ of $Z$ is high (or even infinite)
dimensional and $F$ is a nonlinear function. For instance, the
semilinear partial differential equation of the type $u_t =
\mathcal{A} u + \psi(u,\nabla u,\nabla^2 u, \dots)$ with periodic
boundary conditions in an extended domain, often appearing in pattern
formation dynamics of spatially extended systems, can be written as an
infinite dimensional system of ordinary differential equations (ODEs)
in Fourier space in the form of \eqref{eq:dyn:sys:prototype} in which
case $F$ depends on the operator $\mathcal{A}$ and the function
$\psi$. As we are only interested in the evolution of a few selected
degrees of freedom, i.e.\ only some components of the full dynamics
$Z$ solving \eqref{eq:dyn:sys:prototype}, we assume that one can
separate these resolved degrees of freedom (RDoF) in the dynamical
system from the unresolved degrees of freedom (UDoF). The choice of
RDoF and UDoF is a part of our modeling strategy. Standard examples
include systems with well-separated time scales, e.g. the
decomposition between climate and weather degrees of freedom in
atmosphere-ocean science and the use of reaction coordinates in the
study of chemical kinetics or in molecular dynamics. For such systems,
one decomposes the state space into subspaces $\mathcal{X}$ and
$\mathcal{Y}$ that contain the RDoF and UDoF, respectively:
\begin{equation*}
\mathcal{Z} = \mathcal{X} \oplus \mathcal{Y}\;,
\end{equation*}
with $\mbox{dim}(\mathcal{X}) \ll \mbox{dim}(\mathcal{Y})$
typically. We also introduce the projection operators onto these
spaces $P : \mathcal{Z} \mapsto \mathcal{X}$ and $(I-P): \mathcal{Z}
\mapsto \mathcal{Y}$, respectively. Let now $X$ be the projection of
$Z$ onto the space of $\mathcal{X}$, i.e.\ $X = P Z$. Then we
postulate the existence of a reduced coarse-grained stochastic model
describing the evolution of $X$ alone.
Here we assume that the stochastic model for $X$ is given via a
stochastic differential equation (SDE):
\begin{equation}
  dX = f(X)\,dt + \sqrt{g(X)}\,dW_t\;,\label{eq:reduced:prototype}
\end{equation}
where $W$ denotes a standard Brownian motion of dimension equal to
$\mbox{dim}(\mathcal{X})$. Once the coarse-grained model
\eqref{eq:reduced:prototype} is identified, it can be a used to study
the dynamic characteristic features of the full system
\eqref{eq:dyn:sys:prototype}. Indeed, its low-dimensionality and
simplicity makes it particularly accessible for both rigorous and
computational treatment; see \cite{Pradas2011,Pradas2012,Schmuck2013}
for examples. For many practically relevant cases however, and as we
emphasized earlier, it is not possible to derive a coarse-grained
model \eqref{eq:reduced:prototype} analytically, because of the
complexity of the underlying full system or simply because the full
model \eqref{eq:dyn:sys:prototype} is not completely known.
Consequently, the only way to obtain a coarse-grained model in such a
situation is to use observations, e.g.\ experimental and/or simulation
data, of the full dynamics projected onto the subspace $\mathcal{X}$,
i.e.\ onto of RDoF. That is, it is desirable to identify the
coarse-grained SDE model \eqref{eq:reduced:prototype} in a data-driven
fashion.

An important class of dynamical systems for which coarse-grained
equations of the form \eqref{eq:reduced:prototype} are known to exist,
is when the dynamical system \eqref{eq:dyn:sys:prototype} is given as
system of SDEs with two widely separated time scales. Such systems are
a natural testbed for data-driven coarse-graining techniques, as one
has explicit information about the coarse-grained model. Specifically,
let us consider the following as a prototypical multiscale system
\begin{subequations}
  \begin{align}
    dX^\varepsilon &= \biggl(\frac{1}{\varepsilon}a_0(X^\varepsilon,Y^\varepsilon)+a_1(X^\varepsilon,Y^\varepsilon)\biggr)\,dt + \alpha_0(X^\varepsilon,Y^\varepsilon)\, dU_t + \alpha_1(X^\varepsilon,Y^\varepsilon)\, dV_t\;,\label{eq:sde:generic:fastslow:slow}\\
    dY^\varepsilon &=  \biggl(\frac{1}{\varepsilon^{2}}b_0(X^\varepsilon,Y^\varepsilon) +
    \frac{1}{\varepsilon}b_1(X^\varepsilon,Y^\varepsilon)\biggr)\,dt +
    \frac{1}{\varepsilon}\beta(X^\varepsilon,Y^\varepsilon)\,dV_t\;,\label{eq:sde:generic:fastslow:fast}
  \end{align}
  \label{eq:sde:generic:fastslow}%
\end{subequations}%
with $\varepsilon\ll 1$ controlling the time scale separation. That
is, $X^\varepsilon$ denotes the degrees of freedom we are interested
in (i.e.\ the RDoF) and for which we would like to obtain a
coarse-grained model describing the evolution of $X^\varepsilon$
independent of $Y^\varepsilon$ as $\varepsilon\ll 1$. Mathematically,
the derivation of such coarse-grained models can be made rigorous in
the limit of $\varepsilon\rightarrow 0$ using averaging and
homogenization techniques; see e.g.\ \cite{Pavliotis2008book} and the
references therein for details. In particular, the slow process
$X^\varepsilon$ converges weakly in $C([0,T],\R^d)$ to $X$ solving an
SDE of the form \eqref{eq:reduced:prototype}:
\begin{equation}
  dX = f(X)\,dt + \sqrt{g(X)}\,dW_t\;.\label{eq:sde:generic:effective}
\end{equation}
The drift and diffusion coefficients (i.e.\ the functions $f$ and $g$)
can be formally derived using standard results from homogenization
theory. A data-driven coarse-graining strategy would then be to use
available observations of the multiscale system, specifically of
$X^\varepsilon$ in \eqref{eq:sde:generic:fastslow}, to identify the
coarse-grained model \eqref{eq:sde:generic:effective} by inferring the
functions $f$ and $g$.

Often it is possible to justify proposing a coarse-grained equation
with a particular structure based on theoretical arguments or previous
experience with similar systems. In these cases the inference problem
for $f$ and $g$ in Eq.~\eqref{eq:sde:generic:effective} reduces to
estimating unknown parameters in the SDE. There is a vast and rich
literature on the parametric inference problem for SDEs; see
\cite{PrakasaRao1999,Kutoyants2004,Liptser2010} for instance. For a
data-driven coarse-graining approach for
Eq.~\eqref{eq:sde:generic:effective} based on observations from
Eq.~\eqref{eq:sde:generic:fastslow} it turns out, however, that
commonly used estimators can be biased due to small scale effects in
the observations. In fact, estimators, such as the maximum likelihood
estimator and the quadratic variation of the path estimator, are
highly sensitive to the scale separation. While these estimators do
converge (as $\varepsilon\rightarrow 0$) to the parameters in the
coarse-grained model on the shorter advective time scale, they become
biased on the longer diffusive time scale
\cite{Pavliotis2007,Papavasiliou2009}. The systematic bias due to
multiscale effects on the diffusive time scale can be reduced by
subsampling the data at an appropriate rate. However, the idea of
subsampling does not necessarily lead to an efficient algorithm that
can be used by practitioners, because the optimal sampling rate is
known only for very simple systems (see e.g.\
\cite{Zhang2005,Azencott2010,Azencott2011}) and since, furthermore,
subsampling the data increases the variance of the estimator. A
satisfactory algorithm for fitting a coarse-grained SDE to data based
on the idea of subsampling at the optimal rate combined with an
appropriate variance reduction step has been developed only for some
simple systems used in econometrics \cite{Zhang2005}. To our knowledge
such a methodology has not been developed and implemented for problems
arising in the natural sciences, such as in molecular dynamics or in
statistical physics for example. In addition to the problem of
typically not knowing the optimal subsampling rate, the numerical
experiments in \cite{Pavliotis2007} moreover indicate that the optimal
subsampling rate can vary between parameters in the same
coarse-grained model. Related work that investigates the problem of
parametric inference combined with subsampling techniques in various
settings can be found e.g.\ in
\cite{Cotter2009,Olhede2009,Crommelin2011,Crommelin2012}, while
parametric inference for multiscale problems with vanishing noise is,
e.g., also \cite{Spiliopoulos2013}. Similar consistency questions
arise also in fields other than parametric inference, including
problems in stochastic filtering and stochastic control for SDEs with
multiple scales \cite{Imkeller2013,Zhang2014}.

Related data-driven approaches have also been studied in the context
of numerical methods for SDEs with multiple time scales, i.e.\ for
systems of the form \eqref{eq:sde:generic:fastslow}. We mention in
particular the heterogeneous multiscale method
\cite{Vanden-Eijnden2003,E2005}, which is based on the idea of
evolving the solution of the reduced coarse-grained equation, when the
coefficients in the coarse-grained equation are being evaluated ``on
the fly'' by running short runs of the underlying fast
dynamics. Similar ideas have been proposed in the framework of the
equation-free methodology introduced by Kevrekidis and collaborators
(see e.g.\
\cite{Theodoropoulos2000,Kevrekidis2003,Kevrekidis2004,Kevrekidis2009}),
where a coarse-grained model is evolved using appropriately
initialized simulation on short time scales of the full multiscale
system without knowing the coarse-grained equation in closed form,
making this methodology in principle also applicable for more general
problems, such as kinetic equations. As such, these techniques can be
viewed as a hybrid between numerical analysis and statistical
inference.

To accurately infer coarse-grained models from observations of a
multiscale system, one has to resort to alternative estimation
methodologies, which are robust with respect to the multiscale effects
of the dynamics. The present study is motivated by a recently
introduced estimation methodology which demonstrated how to bypass the
need to subsample data~\cite{Krumscheid2013}.  In the form as proposed
in \cite{Krumscheid2013} this methodology is, however, only applicable
for observations where an ensemble of short trajectories for multiple
initial conditions is available; a design common in many
computer-based simulations. In most real world experiments, such as in
molecular dynamics simulations one typically has access only to a
single long time series. The goal of the present work is therefore to
generalize and appropriately extend the methodology developed in
\cite{Krumscheid2013}, so that it can be used for an observation
design where only one long time series is available. Furthermore, we
demonstrate by means of numerical experiments that the proposed
inference methodology works well for various quite general dynamical
systems of the form \eqref{eq:dyn:sys:prototype}, for which a
coarse-grained model of the form \eqref{eq:reduced:prototype} is known
to exist.

The rest of the paper is organized as follows. In Section
\ref{sec:estimator} we follow the general procedure of
\cite{Krumscheid2013} and present the necessary generalizations and
extensions required for the case of the observational design of a
single time series. Specifically, we will focus on problems where both
the drift function $f$ and the diffusion coefficient $g$ depend
linearly on an unknown parameter vector. This setting covers, for
example, the cases where $f$ and $g$ can be expressed as an
appropriate series (e.g.\ Taylor of Fourier) of known
functions. Moreover, we believe that the ideas developed in this work
for the case of linear parameter dependency can also be instrumental
for the nonlinear case. To demonstrate the effectiveness of the
developed methodology we apply it to a number of selected examples,
which we discuss in Section \ref{sec:numerics}. Specifically, we use
the estimation methodology to identify coarse-grained models for
Brownian motion in a two-scale potential (i.e.\ a stochastic
multiscale systems), for a deterministic system exhibiting chaos, for
a Kac--Zwanzig model, and for a deterministic model for Brownian
motion.  Finally, Section \ref{sec:conclusion} offers a summary and
discussion of our results.

%
%
\section{Estimators for coarse-grained models}
\label{sec:estimator}
We outline here a general methodology that can be used to estimate
parameters in SDEs based on a single trajectory of discrete time
observations.  For the sake of clarity, we first outline the
derivation of the estimator for the case where no multiscale effects
are present. To this end we derive an estimating equation in a
continuous time setting, which will relate the unknown parameters to
statistical properties of the solution to the SDE and discuss how to
obtain parametric estimators from it. To obtain a functional relation
between unknown parameters and statistical properties of the model, in
Section \ref{sec:estimator:esteqn} we follow the methodology outlined
in \cite{Krumscheid2013} and generalize it appropriately. Most of the
examples we are interested in are such that the coarse-grained model
is one-dimensional, see Section \ref{sec:numerics}, i.e.\ we focus on
the case of a scalar diffusion process. It is, however, worthwhile to
remark that the our derivation can be readily extended to the
multidimensional case. Moreover, we discuss modifications and
discretizations to the continuous time estimating equation to account
for observations which are available in the from of a time series
before discussing the coarse-graining scenario.

\subsection{Estimating equation}
\label{sec:estimator:esteqn}
Consider the scalar-valued It{\^o} SDE
\begin{equation}
  dX = f(X)\,dt + \sqrt{g(X)}\,dW_t\;,\quad X(0) = \xi\;,\label{eq:sde:generic}
\end{equation}
on some finite time interval $[0,T]$, $T>0$, with $W$ denoting a
standard one-dimensional Brownian motion. We assume that both the
drift function $f$ and the diffusion function $g$ are such that
Eq.~\eqref{eq:sde:generic} has a unique strong solution on $[0,T]$;
details are given in~\cite{Karatzas1991,Oksendal2003}.  Let us denote
by $X_\xi(t)$ the solution of Eq.~\eqref{eq:sde:generic} at time $t$
started in $\xi$ at time zero, i.e.\ $X_\xi(0) = \xi$. Moreover,
denote by $\mathcal{L}$ the generator associated with
\eqref{eq:sde:generic}, i.e.\ $\mathcal{L} := f\frac{d}{dx} +
\frac{1}{2}g\frac{d^2}{dx^2}$. Then, It{\^o}'s formula together with
the martingale property of the stochastic integral implies that
\begin{equation}
  \E\Bigl(\phi\bigl(X_\xi(t)\bigr)\Bigr) - \phi(\xi)
   = \int_0^t\E\Bigl((\mathcal{L}\phi)\bigl(X_\xi(s)\bigr)\Bigr)\,ds\;,\label{eq:ito:est:phi}
\end{equation}
for any $\phi\in C^2(\R)$ and deterministic initial condition
$\xi$. For the sake of completeness, we remark that other commonly
used notations for $\E\bigl(\phi\bigl(X_\xi(t)\bigr)\bigr)$ are
$\E\bigl(\phi\bigl(X(t)\bigr)\bigl.\bigr\vert X(0) = \xi\bigr)$ or
$\E_\xi\bigl(\phi\bigl(X(t)\bigr)\bigr)$, and that
Eq.~\eqref{eq:ito:est:phi} is also known as Dynkin's formula
\cite[Ch.\ $7.4$]{Oksendal2003}.

In this work we follow a semiparametric approach for the
parametrization of Eq.~\eqref{eq:sde:generic}. That is, we assume that
both $f$ and $g$ depend on an unknown parameter vector $\theta\equiv
(\theta_1,\dots,\theta_n)^T\in\R^n$, $n\in\N$, which we wish to
determine from observations. Specifically, we consider
\begin{equation}
  f(x)\equiv f(x;\theta) := \sum_{j=1}^n\theta_jf_j(x)\quad\textrm{and}\quad
  g(x)\equiv g(x;\theta) := \sum_{j=1}^n\theta_jg_j(x)\;,\label{eq:parameterization}
\end{equation}
with some known functions $f_j$ and $g_j$, $1\le j\le n$. That is,
both $f$ and $g$ can depend on the same parameter. If this not the
case however, one can think of the first $k$, say, components of the
vector $\theta$ parametrizing the drift function $f$ while the
remaining $n-k$ components the diffusion function $g$, and setting
$f_j=0$ for $k<j\le n$ as well as $g_j=0$ for $1\le j\le k$. For the
numerical examples in Section \ref{sec:numerics} we will have that $f$
and $g$ are polynomials of some degree, so that $f_j$ and $g_j$ will
be appropriate monomials, respectively. After substituting
\eqref{eq:parameterization} into \eqref{eq:ito:est:phi} and
rearranging terms, we arrive at
\begin{equation}
  \E\Bigl(\phi\bigl(X_\xi(t)\bigr)\Bigr) - \phi(\xi)
  = \sum_{j=1}^n\theta_j\int_0^t\E\Bigl((\mathcal{L}_j\phi)\bigl(X_\xi(s)\bigr)\Bigr)\,ds\;,
   \label{eq:ito:est:phi:param}
\end{equation}
with $\mathcal{L}_j := f_j\frac{d}{dx} +
\frac{1}{2}g_j\frac{d^2}{dx^2}$. To write this estimating equation
\eqref{eq:ito:est:phi:param} in a more compact manner, we define the
following component functions for any fixed time $t\in[0,T]$ and any
fixed function $\phi$,
\begin{equation*}
  b_c(\xi) :=  \E\Bigl(\phi\bigl(X_\xi(t)\bigr)\Bigr) - \phi(\xi) \in\R\quad\text{and}\quad
  a_c(\xi) := \Bigl(\int_0^t\E\Bigl((\mathcal{L}_j\phi)\bigl(X_\xi(s)\bigr)\Bigr)\,ds\Bigr)_{1\le j\le n}\in\R^n\;,
\end{equation*}
which highlight the dependency on $\xi$. Using these definitions,
Eq.~\eqref{eq:ito:est:phi:param} reduces to
\begin{equation}
a_c(\xi)^T\theta = b_c(\xi)\;,\label{eq:fun:form:param:single}
\end{equation}
which is underdetermined for $n>1$.  To make this identity useful
nonetheless, we exploit the fact that
Eq.~\eqref{eq:fun:form:param:single} is valid for any $\xi$; a
technique that has already been used successfully in
\cite{Krumscheid2013}. We now introduce the concept of \textit{trial
  points}: as we work in an observation framework where only one time
series is available, we denote by $\xi$ the trial point instead of
initial condition to avoid confusion with the initial condition of the
time series; see also
Section~\ref{sec:estimator:esteqn:mod:condexp}. By considering a
finite sequence of trial points ${(\xi_{i})}_{1\le i\le m}$, we can
assemble a system of linear equations, solved by the parameter vector
$\theta$:
\begin{equation}
  A\theta = b\;,\label{eq:fun:form:param:sys}
\end{equation}
where $A := \bigl(a_c(\xi_{i})^T\bigr)_{1\le i\le m}\in\R^{m\times n}$
and right-hand side $b:= \bigl(b_c(\xi_{i})\bigr)_{1\le i\le
  m}\in\R^{m}$. Since this linear system does not have a unique
solution in general, we define the estimator of $\theta$ based on $A$
and $b$ as the least squares solution of $A\theta = b$ with minimum
norm:
\begin{equation}
  \hat{\theta} := \argmin_{x\in\mathcal{S}}\norm{x}_2^2\;,
  \quad\mathcal{S}:=\bigl\{x\in\R^n\colon \norm{Ax-b}_2^2 = \min\bigr\}\;.\label{eq:param:lsp}
\end{equation}
At this point, we can still exploit the degree of freedom for choosing
$\phi$ in Eq.~\eqref{eq:ito:est:phi:param} freely. Motivated by
\cite{Krumscheid2013} where approximations of the first and second
moment provided very accurate estimates of $\theta$, we use
$\phi(x):=x+x^2$ throughout this work. In fact, the two-step
estimation approach for $\theta$ presented in the aforementioned work
can be recovered as a special case of the procedure outlined
here. Indeed, using $\phi(x) = x$ causes
Eq.~\eqref{eq:ito:est:phi:param} to degenerate to an equation not
containing any parameters characterizing the diffusion function
$g$. This then yields an estimator for the drift parameters
only. After this first step, we substitute the obtained estimators
into the parametrization of $f$. Repeating then the same steps with
the function $\phi(x) = x^2$ gives an estimator of the remaining
parameters determining $g$ and the two-step scheme is completed.
Finally, we mention that other choices for the function $\phi$ may
work as well to construct the estimator \eqref{eq:param:lsp}.  In view
of the preceding discussion it is clear however, that $\phi$ cannot be
chosen arbitrarily as it has to be such that
Eq.~\eqref{eq:ito:est:phi:param} still depends on all unknown
parameters that we want to estimate. The used function $\phi(x)=x+x^2$
thus appears to be the obvious candidate for various problems, while
other choices of $\phi$ may be problem dependent. A more systematic
study of how to choose $\phi$ in the context of the rigorous
convergence analysis will be presented in \cite{Krumscheid_pre}.

\subsection{Modifications due to discrete time observations}
\label{sec:estimator:esteqn:mod} Recall that we seek to determine an
approximation of the parameter vector $\theta$ in
Eq.~\eqref{eq:sde:generic} with parametrization
\eqref{eq:parameterization}, based on a trajectory of discrete time
observations. That is, we have access to $N$ data $\mathbb{X}_N :=
{\bigl(X(t_k)\bigr)}_{1\le k \le N }$ with $t_k = (k-1)h$, where
$h=T/(N-1)$. A constant sampling rate $h$ is assumed here merely for
simplicity and the proposed methodology can be readily extended to the
case of non-constant sampling rates.  To apply the methodology outlined
above, we have to carry out two essential modifications to the purely
continuous framework \eqref{eq:ito:est:phi:param}. Firstly, we have to
estimate the conditional expectations of the form
$\E\bigl(\varphi\bigl(X_\xi(\tau)\bigr)\bigr)$ based on
$\mathbb{X}_N$. Secondly, we have to replace the temporal integrals
with discrete versions. A detailed algorithmic description of the
estimation procedure for discrete time observations based on these
modifications is presented in Section
\ref{sec:estimator:esteqn:mod:algo}.

\subsubsection{Estimating the conditional expectation}
\label{sec:estimator:esteqn:mod:condexp}
Throughout the estimation procedure, we have to approximate
conditional expectations of the form
$\E\bigl(\varphi\bigl(X_\xi(\tau)\bigr)\bigr)$ for multiple values of
the trial point $\xi$. The available time series $\mathbb{X}_N$
provides, however, only one initial condition which we cannot
influence nor manipulate; thus the necessity to distinguish between
trial point and initial condition. A way out of this predicament is
possible when the time series (i.e.\ the discrete time process) is
stationary and sufficiently mixing so that
\begin{equation*}
  \cov\bigl(X(t),X(t+kh)\bigr) \le C \rho^k\;,
\end{equation*}
for some finite $C>0$ and $\rho\in[0,1[$, which we will assume from
now on; see e.g.~\cite{Bosq1998,Fan2003} for further details in the
inference context and refer to e.g.~\cite{Doukhan1994,ChenX2010} and
the references therein for a discussion of sufficient conditions on
the drift function and the diffusion coefficient in SDE model
\eqref{eq:sde:generic}. Related conditions on the covariance as a
function of the lag $k$ have also been used in other works on
parametric inference for diffusion processes; see \cite{Azencott2011}
for instance. Intuition in this case then suggests to sequentially
search the time series $\mathbb{X}_N$ for the value of the trial point
$\xi$ and then to approximate the expectation by averaging over the
events $\varphi(X)$ at $\tau$ time units after the occurrences of
$\xi$ in $\mathbb{X}_N$. A technique which makes this approximation
idea precise is the class of so-called local polynomial kernel
regression estimators \cite{Fan2003}. Recall that the sampling time of
the time series $\mathbb{X}_N$ is $h$. For a shift by $\tau >0$ time
units to be well-defined, we require that $\tau = lh$, for some
$l\in\{1,2,\dots,N-1\}$ and for such a $\tau$ we set $N_\tau =
N-\tau/h\in\N$. Then the simplest regression estimator (locally
constant) yields the approximation
\begin{equation}
   {\Bigl.\E\Bigl(\varphi\bigl(X_\xi(\tau)\bigr)\Bigr)\Bigr\vert}_{\tau=lh}\approx
   \frac{\sum_{k=1}^{N_\tau}\varphi\bigl(X(t_{k+l})\bigr)K\Bigl(\frac{X(t_k)-\xi}{\kappa_{N_\tau}}\Bigr)}
   {\sum_{k=1}^{N_\tau}K\Bigr(\frac{X(t_k)-\xi}{\kappa_{N_\tau}}\Bigr)}\;,
 \label{eq:approx:exp:mixing}
\end{equation}
which is also known as the Nadaraya--Watson estimator
\cite{Nadaraya1964,Watson1964}. Here $K$ is an appropriately chosen
kernel, and $0<\kappa_{N_\tau}$ denotes the bandwidth which decays to
zero as $N_\tau\rightarrow 0$ at a rate depending on the sense of
convergence in Eq.~\eqref{eq:approx:exp:mixing}; details are given
in~\cite{Bosq1998}.  Throughout this study we select the Gaussian
kernel $K(x) := \exp{(-x^2/2)}/\sqrt{2\pi}$ for convenience, but we
remark that other choices are possible.

Upon defining $w_{N_\tau,k}(\xi) :=
K\bigl((X(t_k)-\xi)/\kappa_{N_\tau}\bigr)/
\sum_{k=1}^{N_\tau}K\bigr((X(t_k)-\xi)/\kappa_{N_\tau}\bigr)$, one can
rewrite the regression estimator, i.e.\ the right-hand side in
Eq.~\eqref{eq:approx:exp:mixing}, as
$\sum_{k=1}^{N_\tau}w_{N_\tau,i}(\xi)\varphi\bigl(X(t_{k+l})\bigr)$. That
is, the regression estimator is given as a weighted average with
non-identical weights $w_{N_\tau,k}(\xi)$. Let us finally note that if
the trial point $\xi$ is such that the denominator of the regression
estimator in Eq.~\eqref{eq:approx:exp:mixing} is zero (roughly
speaking this happens if $\xi$ is not in the support of the stationary
density of $\mathbb{X}_N$), then we set $w_{N_\tau,k}(\xi) = 1/N_\tau$
instead for well-posedness (see also Sect.\
\ref{sec:estimator:esteqn:mod:algo} below). However, one should ensure
that this event is avoided by selecting the trial points
appropriately, otherwise the estimator's approximation accuracy would
deteriorate due to incorporating unfeasible information. As the
regression estimator in Eq.\ \eqref{eq:approx:exp:mixing} essentially
averages over the events $\varphi\bigl(X(t_{k+l})\bigr)$ for which
$X(t_k)\approx \xi$, one should moreover try to ensure that the trial
point $\xi$ is located in a region where most of the observations are
located in order to average over a sufficiently large sample; see
Section \ref{sec:numerics} for a detailed description of how to chose
the trial points in practice.

\subsubsection{Temporal integrals}
The integrands of the temporal integrals in
Eq.~\eqref{eq:ito:est:phi:param} are precisely the conditional
expectations discussed above. Let $u(\tau) :=
\E\bigl(\varphi\bigl(X_\xi(\tau)\bigr)\bigr)$ be such an expectation
for a fixed trial point $\xi$ and function $\varphi$. To replace the
temporal integral of $u$ over $[0,t]$ by a discrete version in
\eqref{eq:ito:est:phi:param}, we use the composite trapezoidal rule
with $n_h$ equally spaced ($n_h = t/h$) subdivisions:
\begin{equation}
  \int_0^t u(s)\,ds \approx\frac{h}{2} \biggl(u(0) + u(t) +  2\sum_{l=1}^{n_h-1}u(lh) \biggr)\;.\label{eq:trapz:rule}
\end{equation}
The choice of an equally spaced subdivision of $[0,t]$ where the
division length coincides with the sampling rate $h$ of the available
time series $\mathbb{X}_N$ is made for reasons of a consistent
discretization. In fact, it ensures that the time points $\tau$, say,
at which the integrand $u$ is evaluated, is an integer multiple of
$h$, so that the shifts by $\tau$ time units in the regression
estimator \eqref{eq:approx:exp:mixing} are well-defined. Other
(possibly non-equally spaced) time discretizations, which are
consistent in the sense that each trapezoidal node is an integer
multiple of $h$, are of course possible. Finally, we mention that the
use of trapezoidal rule \eqref{eq:trapz:rule} is motivated by the fact
that the integrands $u$ are replaced by the regression estimators
\eqref{eq:approx:exp:mixing} in practice, for which we cannot expect
to provide sufficient smoothness. Under these conditions the
trapezoidal rule is advantageous over higher order methods since
higher order derivatives, as used in classical Taylor expansion based
arguments, are not continuous \cite{Cruz-Uribe2002}.

\subsubsection{An algorithmic description for discrete time
  observations}
\label{sec:estimator:esteqn:mod:algo} To illustrate how the
combination of these approximations can be used to apply the developed
methodology to discrete time observations, we present a detailed
pseudocode in Algorithm \ref{alg:est}.
\begin{algorithm}[t]
\begin{algorithmic}[1]
\REQUIRE $0<t$ such that $t/h\in\N$, $0<h$, $\Xi\in\R^m$, and $\mathbb{X}_N\in\R^N$
\STATE $l \leftarrow \frac{t}{h}$
\FOR{$i=1$ to $m$}
\STATE $\xi \leftarrow \Xi_i$
\FOR{$j=1$ to $n$}
\STATE $u_{j,0} \leftarrow f_j(\xi)(1+2\xi) + g_j(\xi)$
\ENDFOR
\FOR{$k=1$ to $l$}
\STATE $X \leftarrow \mathbb{X}_N(1:N-k)$
\STATE $Y \leftarrow \mathbb{X}_N(1+k:N)$
\FOR{$j=1$ to $n$}
\STATE $u_{j,k} \leftarrow \texttt{nwe}\bigl(X,f_j(Y)(1+2Y) + g_j(Y),\xi\bigr)$
\ENDFOR
\ENDFOR
\STATE $X \leftarrow \mathbb{X}_N(1:N-l)$
\STATE $Y \leftarrow \mathbb{X}_N(1+l:N)$
\STATE $b_i \leftarrow \texttt{nwe}\bigl(X, Y+Y^2 ,\xi\bigr) - (\xi+\xi^2)$
\FOR{$j=1$ to $n$}
\STATE $A_{i,j} \leftarrow \frac{h}{2}\bigl(u_{j,0} + u_{j,l} + 2\sum_{k=1}^{l-1} u_{j,k}\bigr)$
\ENDFOR
\ENDFOR
\STATE $\theta \leftarrow A^{+}b$
\RETURN $\theta$
\end{algorithmic}
\caption{Algorithmic description of the introduced estimation
  procedure.}
\label{alg:est}
\end{algorithm}
Here we assume that a parametrization for both drift function and
diffusion function has been fixed by choosing $f_j$ and $g_j$ in
Eq.~\eqref{eq:parameterization}, for $1\le j\le n$. The input
arguments of Algorithm \ref{alg:est} are the time series
$\mathbb{X}_N$ of $N$ discrete time observations corresponding to a
constant sampling rate $h$, the $m$ trial points $\Xi$, and the time
$t$ controlling the temporal integration in
\eqref{eq:ito:est:phi:param}, which is assumed to be an integer
multiple of $h$ (cf.\ Section
\ref{sec:estimator:esteqn:mod:condexp}). We note that we use the colon
notation \cite[Ch.\ $1.1.8$]{Golub1996} in lines $8,9$ and $14,15$ to
select several components of a vector at once, so that we can suppress
additional iteration details. Similarly, the application of a function
defined on $\R$ to a vector (such as in lines $11$ and $16$) is
understood componentwise. We emphasize that the statement $\theta
\leftarrow A^{+}b$ in line $21$ is merely meant as a formal notation
for computing the least squares solution of $A\theta = b$ with minimum
norm. In fact, in this work we use a QR factorization with column
pivoting to solve the least squares problem but other choices are
possible, typically depending on the rank of $A$; see, e.g.,
\cite[Ch.\ $5$]{Golub1996}.  Furthermore, the procedure \texttt{nwe}
(called in lines $11$ and $16$) implements the Nadaraya--Watson
estimator \eqref{eq:approx:exp:mixing} for the approximation of
conditional expectations and its detailed pseudocode is given in
Algorithm \ref{alg:est:nwe}.
\begin{algorithm}[t]
\begin{algorithmic}[1]
\REQUIRE $X,Y\in\R^M$ and $\xi\in\R$
\STATE $\kappa \leftarrow \argmin_{\delta>0} \Bigl(\frac{1}{\delta M^2 \sqrt{2}}\sum_{i,j=1}^M K\bigl(\frac{X_i-X_j}{\delta\sqrt{2}}\bigr) - \frac{2}{M(M-1)}\sum_{i=1}^M\sum_{i\not= j} K\bigl(\frac{X_i-X_j}{\delta}\bigr)\Bigr)$
\IF{$\sum_{i=1}^M K\bigl(\frac{X_i-\xi}{\kappa}\bigr) = 0$}
\STATE $u \leftarrow \frac{1}{M}\sum_{i=1}^M Y_i$
\ELSE
\STATE $u \leftarrow \frac{\sum_{i=1}^M Y_i K\bigl(\frac{X_i-\xi}{\kappa}\bigr)}{\sum_{i=1}^M K\bigl(\frac{X_i-\xi}{\kappa}\bigr)}$
\ENDIF
\RETURN $u$
\end{algorithmic}
\caption{Pseudocode of the \texttt{nwe} procedure used in Algorithm
  \ref{alg:est} to approximate conditional expectations via the
  Nadaraya--Watson estimator.}
\label{alg:est:nwe}
\end{algorithm}
Its input arguments are two lists $X$, $Y$ of the same length as well
as the trial point $\xi$ and the algorithm returns an approximation of
$\E(Y\vert X=\xi)$. In the pseudocode presented here, we use the least
squares cross validation for a data-driven bandwidth selection (line
$1$). We mention, however, that this selection technique is used here
merely for the sake of a compact notation and several other methods
can be used alternatively \cite[Ch.\ $8.5$]{Kroese2011}. For the
numerical examples discussed in Section \ref{sec:numerics} we tried
different bandwidth selection methods (not shown) but did not observe
any significant differences. We also note that there exist efficient
computational strategies to evaluate the term in brackets in line $1$
of Algorithm \ref{alg:est:nwe} via fast Fourier transform related
approaches.

\subsection{Estimators for coarse-grained models of multiscale systems}
\label{sec:estimator:multiscale}
A central goal of this study is to identify a coarse-grained model
based on observations of a multiscale system. Specifically, we
consider the prototypical multiscale system
\eqref{eq:sde:generic:fastslow}, i.e.\
\begin{subequations}
  \begin{align}
    dX^\varepsilon &= \biggl(\frac{1}{\varepsilon}a_0(X^\varepsilon,Y^\varepsilon)+a_1(X^\varepsilon,Y^\varepsilon)\biggr)\,dt + \alpha_0(X^\varepsilon,Y^\varepsilon)\, dU_t + \alpha_1(X^\varepsilon,Y^\varepsilon)\, dV_t\;,\label{eq:sde:generic:fastslow2:slow}\\
    dY^\varepsilon &=  \biggl(\frac{1}{\varepsilon^{2}}b_0(X^\varepsilon,Y^\varepsilon) +
    \frac{1}{\varepsilon}b_1(X^\varepsilon,Y^\varepsilon)\biggr)\,dt +
    \frac{1}{\varepsilon}\beta(X^\varepsilon,Y^\varepsilon)\,dV_t\;,\label{eq:sde:generic:fastslow2:fast}
  \end{align}
  \label{eq:sde:generic:fastslow2}%
\end{subequations}%
equipped with appropriate initial conditions on the time interval
$[0,T]$, where $U$ and $V$ denote independent Brownian motions, and
$\varepsilon>0$ is a small parameter controlling the scale
separation. Here we assume that $1=\dim{(\mathcal{X})}$, while
$\dim{(\mathcal{Y})}$ is arbitrary, so that the coarse-grained model
\begin{equation}
  dX = f(X)\,dt + \sqrt{g(X)}\,dW_t\;,\label{eq:sde:generic:effective2}
\end{equation}
is also one-dimensional; $W$ is a standard one-dimensional Brownian
motion. In fact, using results from homogenization theory one can
rigorously show that the process $X^\varepsilon$ solving
\eqref{eq:sde:generic:fastslow2:slow} converges weakly in
$C([0,T],\R)$ to the process $X$ solving
\eqref{eq:sde:generic:effective2} as $\varepsilon\rightarrow 0$,
provided that the fast process $Y^\varepsilon$ is ergodic and the
centering condition is satisfied; see, e.g., \cite{Pavliotis2008book}
and the references therein.

Our data-driven coarse-graining strategy is to use the available
observations of $X^\varepsilon$ solving
\eqref{eq:sde:generic:fastslow2:slow} with $\varepsilon>0$ and
estimate both $f$ and $g$ in \eqref{eq:sde:generic:effective2} using
exactly the same estimation methodology as presented in Section
\ref{sec:estimator:esteqn}. We emphasize that we are facing a problem
of model misspecification now: fitting model
\eqref{eq:sde:generic:effective2} to observation from
\eqref{eq:sde:generic:fastslow2:slow} which is not consistent with
model \eqref{eq:sde:generic:effective2}. Parametric inference for
misspecified models in the absence of multiscale effects has been
studied, e.g.\ in \cite[Ch.\ $2.6$]{Kutoyants2004}. Here we expect
that in the limit of infinite scale separation $\varepsilon\rightarrow
0$, the error due to the model misspecification vanishes
\cite{Krumscheid_pre}. Finally, our motivation to resort to the
estimation methodology of Section \ref{sec:estimator:esteqn} also for
this setting stems from recent results in \cite{Krumscheid2013},
where, as described in Section \ref{sec:intro}, a related scheme
demonstrated to be able to accurately estimate the coarse-grained
model from observations of the multiscale system. This favorable
property agrees with our intuition that the estimated model should be
close to the coarse-grained model if the model misspecification is
small, i.e.\ if $\varepsilon\ll 1$.  In fact, in view of the
theoretical results presented in \cite{Krumscheid_pre}, it is expected
that, in the absence of all other error contributions such as the
finite sample size, the estimators converge in a probabilistic sense
(in fact, the convergence is almost surely) to the parameters in the
coarse-grained model in the limit $\varepsilon \rightarrow 0$.

Specifically, in this multiscale setting we have access to $N$
discrete time observations of
Eq.~\eqref{eq:sde:generic:fastslow2:slow}, that is
$\mathbb{X}_N^\varepsilon := {\bigl(X^\varepsilon(t_k)\bigr)}_{1\le k
  \le N }$ with $t_k = (k-1)h$, where $h=T/(N-1)$. Based on the
semiparametric parametrization \eqref{eq:parameterization} of the
drift and the diffusion coefficients in the the coarse-grained model
\eqref{eq:sde:generic:effective2}, the multiscale time series
$\mathbb{X}_N^\varepsilon$ is used to assemble the corresponding
matrix $A^\varepsilon$ and right-hand side $b^\varepsilon$ in
Eq.~\eqref{eq:fun:form:param:sys}. The estimated parameter vector of
the coarse-grained model based on the multiscale data
$\mathbb{X}_N^\varepsilon$ is then given by the least squares solution
of $A^\varepsilon \theta = b^\varepsilon$ with minimum norm. We denote
the estimated parameter vector by $\hat{\theta}^\varepsilon$ to
emphasize the dependency on the multiscale observations
$\mathbb{X}_N^\varepsilon$. It is also worth mentioning that the
estimation procedure is solely derived from the coarse-grained model
\eqref{eq:sde:generic:effective2}. That is, the procedure does not
incorporate any knowledge of the corresponding multiscale system
\eqref{eq:sde:generic:fastslow2}. In addition, we are not assuming
knowledge of the scale separation parameter $\varepsilon$. In other
words, one can view the available time series
$\mathbb{X}_N^\varepsilon$ as obtained purely from a ``black box''
model, which is close to the coarse-grained model
\eqref{eq:sde:generic:effective2} provided that $\varepsilon\ll 1$.

%
%
\section{Numerical experiments}
\label{sec:numerics}

In this section we apply the proposed estimation procedure to several
examples. We focus here on the inference problem for coarse-grained
models based on multiscale observations, for which classical
estimators are expected to fail. In Section\
\ref{sec:numerics:langevin1d} we first investigate a stochastic
multiscale system, namely Brownian motion in a two-scale
potential. The remaining examples are deterministic multiscale
systems, for which we seek to identify a coarse-grained stochastic
model from a single time series. Specifically, we estimate parameters
in the coarse-grained model for a deterministic system exhibiting fast
temporal chaos (Section\ \ref{sec:numerics:langevin1d}), in a
low-dimensional approximation of a large Hamiltonian system (Section\
\ref{sec:numerics:kaczwanzig}), and in an Ornstein--Uhlenbeck process
constructed in a purely deterministic setting (Section\
\ref{sec:numerics:detbm}).  To verify the accuracy of the estimated
parameters in the coarse-grained models, we compare the obtained
estimates with theoretically available ones. We will also address the
question of selecting the time $t$ controlling the temporal integrals
in \eqref{eq:ito:est:phi:param}, in order to uniquely define the
estimation procedure of Section \ref{sec:estimator}. To emphasize the
dependency of the estimated parameter vector based on multiscale data
$\hat{\theta}^\varepsilon$ also on $t$, we use
$\hat{\theta}^\varepsilon\equiv \hat{\theta}_t^\varepsilon$. To
assemble the linear system \eqref{eq:fun:form:param:sys} we use $m=54$
trial point for all examples. As mentioned above already, these points
$\xi$ are selected a-priori for each example such that they cover most
of the range of the time series $\mathbb{X}_N^\varepsilon$.  One
simple and ad hoc way is to identify the region for $\xi$ by defining
$a_N := (1-\nu) \min(\mathbb{X}_N^\varepsilon) +
\nu\max(\mathbb{X}_N^\varepsilon)$ and $b_N :=
\nu\min(\mathbb{X}_N^\varepsilon) + (1 -
\nu)\max(\mathbb{X}_N^\varepsilon)$, for $0<\nu<1/2$. Furthermore, let
$\eta_1, \eta_2, \dots, \eta_m$ be an independent and identically
distributed sequence of random variables following a standard normal
distribution. Then we set $l_m := \min_{1\le i\le m}(\eta_i)$ as well
as $r_m := \max_{1\le i\le m}(\eta_i)$ and select the trial points by
linearly mapping $\eta_i$ to the region of interest $[a_N,b_N]$:
\begin{equation*}
  \xi_i := \frac{a_N-b_N}{l_m-r_m}\eta_i + \frac{l_mb_N-r_ma_N}{l_m-r_m}\;,
\end{equation*}
for $1\le i\le m$, which are then fixed throughout the numerical
experiment. This rather naive procedure worked well for the examples
that follow where we used $\nu = 0.2$, because the trial points $\xi$
are located in regions where most of the observations are, so that
estimates of the conditional expectations are expected to be
accurate. We remark however, that other approaches are possible as
well. In fact, a more systematic way to construct the trial points is
via a resampling method, which exploits the shape of the empirical
distribution function of $\mathbb{X}_N^\varepsilon$. In this case let
$\eta_1, \eta_2, \dots, \eta_m$ be an independent and identically
distributed sequence of random variables following a uniform
distribution on $[0,1]$. Moreover, denote by
$X_{1:N}^\varepsilon,X_{2:N}^\varepsilon,\dots,X_{N:N}^\varepsilon$
the order statistic of $\mathbb{X}_N^\varepsilon$ so that
$X_{1:N}^\varepsilon < X_{2:N}^\varepsilon < \cdots
<X_{N:N}^\varepsilon$. Then the trail points may be constructed via
\begin{equation*}
  \xi_i := X_{k_i:N}^\varepsilon\;,\quad k_i:= \min{\{k\in\N\colon \eta_iN\le k\}}\;,
\end{equation*}
for $1\le i\le m$. Although we did not observe significant differences
between these two trial point selection strategies for the examples
considered here, the second systematic strategy may be advantageous in
more general cases as it is purely data-driven. Finally, we note that
we set the estimation procedure's defining parameters $m$, $T$, and
$h$ in such a way that their error contribution due to approximations
are negligible compared to the scale separation $\varepsilon$. This is
done to focus solely on the estimator's performance under the presence
of multiscale effects in the observations, a scenario where other
estimation techniques fail to be consistent.

\subsection{Brownian particle in a two-scale potential}
\label{sec:numerics:langevin1d}
Let us begin with an example borrowed from \cite{Pavliotis2007}, which
was originally used to investigate the failure of classical parametric
estimation techniques for multiscale diffusion
processes. Specifically, we consider
\begin{equation*}
  dX^\varepsilon = -\frac{d}{dx} V\biggl(X^\varepsilon,\frac{X^\varepsilon}{\varepsilon}\biggr)\,dt + \sqrt{2\sigma}\,dW_t\;,
\end{equation*}
which models the position of a Brownian particle moving in a two-scale
potential $V$ and being affected by thermal noise. Here $W$ denotes a
standard one-dimensional Brownian motion. We investigate the situation when
the two-potential $V$ is given by a large scale part $V_\alpha$ superimposed
with a periodically fluctuating part $p$: $V(x,y) = V_\alpha(x) + p(y)$.
Under this assumption, the multiscale SDE can be written as
\begin{equation}
  dX^\varepsilon = -\biggl(V_\alpha'(X^\varepsilon) + \frac{1}{\varepsilon} p'\biggl(\frac{X^\varepsilon}{\varepsilon}\biggr)\biggr)\,dt + \sqrt{2\sigma}\,dW_t\;. \label{eq:eps:sde:langevin1d}
\end{equation}
Notice that the SDE \eqref{eq:eps:sde:langevin1d} can be rewritten as
a fast/slow system of the form \eqref{eq:sde:generic:fastslow2} by
introducing the auxiliary variable $Y^\varepsilon :=
X^\varepsilon/\varepsilon$.

We consider the case where the fluctuating part $p$ is a smooth
periodic function with period $L$ and the large scale part is a
quadratic potential, i.e.\ $V_\alpha(x) = \alpha x^2/2$. Then, as
$\varepsilon\rightarrow 0$, $X^\varepsilon$ solving
\eqref{eq:eps:sde:langevin1d} converges weakly in $C([0,T], \R)$ to
the solution of the coarse-grained equation
\begin{equation}
  dX = -A X\,dt + \sqrt{2\Sigma}\,dW_t\;,\label{eq:eps:sde:langevin1d:coarse}
\end{equation}
where $A=\alpha L^2/(Z_{+}Z_{-})$ and $\Sigma=\sigma
L^2/(Z_{+}Z_{-})$, with $Z_{\pm} = \int_0^Le^{\pm p(y)/\sigma}\,dy$;
see \cite{Pavliotis2007} for details. We set the fluctuating part to
be $p(y) = \cos{(y)}$. Then the constants $Z_{\pm}$ can be easily
computed so that $A= \alpha/ {I_0(\sigma^{-1})}^{2}$ and $\Sigma =
\sigma/ {I_0(\sigma^{-1})}^{2}$, where $I_0$ denotes the modified
Bessel function of first kind. We note that both parameters in
\eqref{eq:eps:sde:langevin1d:coarse} depend non-trivially on
$\sigma$.

\begin{figure}[t]
  \centering
  \includegraphics[width=0.465\textwidth]{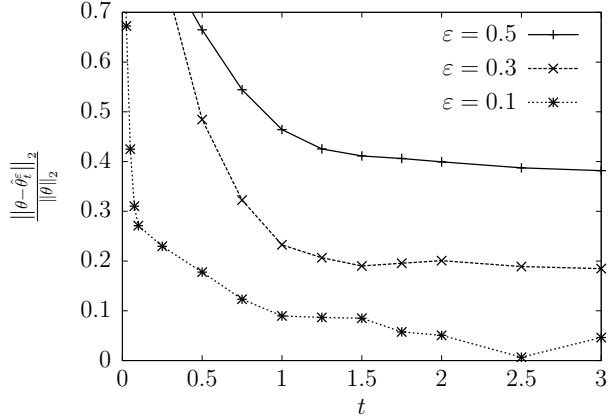}
  \caption[]{Relative error of the estimated parameter vector
    $\hat{\theta}_{t}^\varepsilon$ for
    \eqref{eq:eps:sde:langevin1d:coarse} based on observations of
    \eqref{eq:eps:sde:langevin1d} with $\alpha = 2$, $\sigma=1$, and
    $\varepsilon\in\{0.1,0.3,0.5\}$.}
  \label{figure:langevin}
\end{figure}

To estimate the $n=2$ parameters in
\eqref{eq:eps:sde:langevin1d:coarse}, we choose the functions in the
drift and diffusion parametrization \eqref{eq:parameterization} as
$f_1(x) = x$, $f_2(x) = 0 = g_1(x)$, and $g_2(x) = 2$, with true
parameter vector $\theta = (-A,\Sigma)^T$. The estimate of $\theta$ is
then based on a time series on $[0,1000]$ of the multiscale system
\eqref{eq:eps:sde:langevin1d} with $\alpha = 2$, $\sigma = 1$ for each
$\varepsilon\in\{0.1,0.3,0.5\}$. The time series were obtained by
numerically integrating \eqref{eq:eps:sde:langevin1d} via the
Euler--Maruyama method with step size $h=0.001$ and initial condition
$X^\varepsilon(0) = 0$.  Fig.~\ref{figure:langevin}
shows the relative error of the estimated parameter vector
$\hat{\theta}_{t}^\varepsilon$ as a function of $t$, for each value of
$\varepsilon$ respectively. For all three of them one observes that
while very small values of $t$ result in a large relative error,
increasing $t$ reduces the error significantly. In fact, for $t=1$ and
$\varepsilon = 0.1$ we find a relative error of $10\%$ and for even
larger values of $t$ the relative error drops further significantly
below $5\%$. For $t\ge 3$ (not shown here) the relative error starts
fluctuating around $4\%$ due to discretization errors but remains of
the same order as the scale separation parameter $\varepsilon$, which
is in agreement with the results presented in
\cite{Krumscheid_pre}. Consequently, it is possible to obtain very
accurate estimates of the parameters in the coarse-grained model
\eqref{eq:eps:sde:langevin1d:coarse} based on observations of the
multiscale system \eqref{eq:eps:sde:langevin1d}, once $t$ is
sufficiently large. In fact, by comparing the resulting relative
errors for different values of $\varepsilon$,
Fig.~\ref{figure:langevin} suggests to choose $t$ of $\mathcal{O}(1)$
for a relative error of $\mathcal{O}(\varepsilon)$.

\begin{figure}[t]
  \centering
  \begin{subfigure}[b]{0.465\textwidth}
    \centering
    \includegraphics[width=\textwidth]{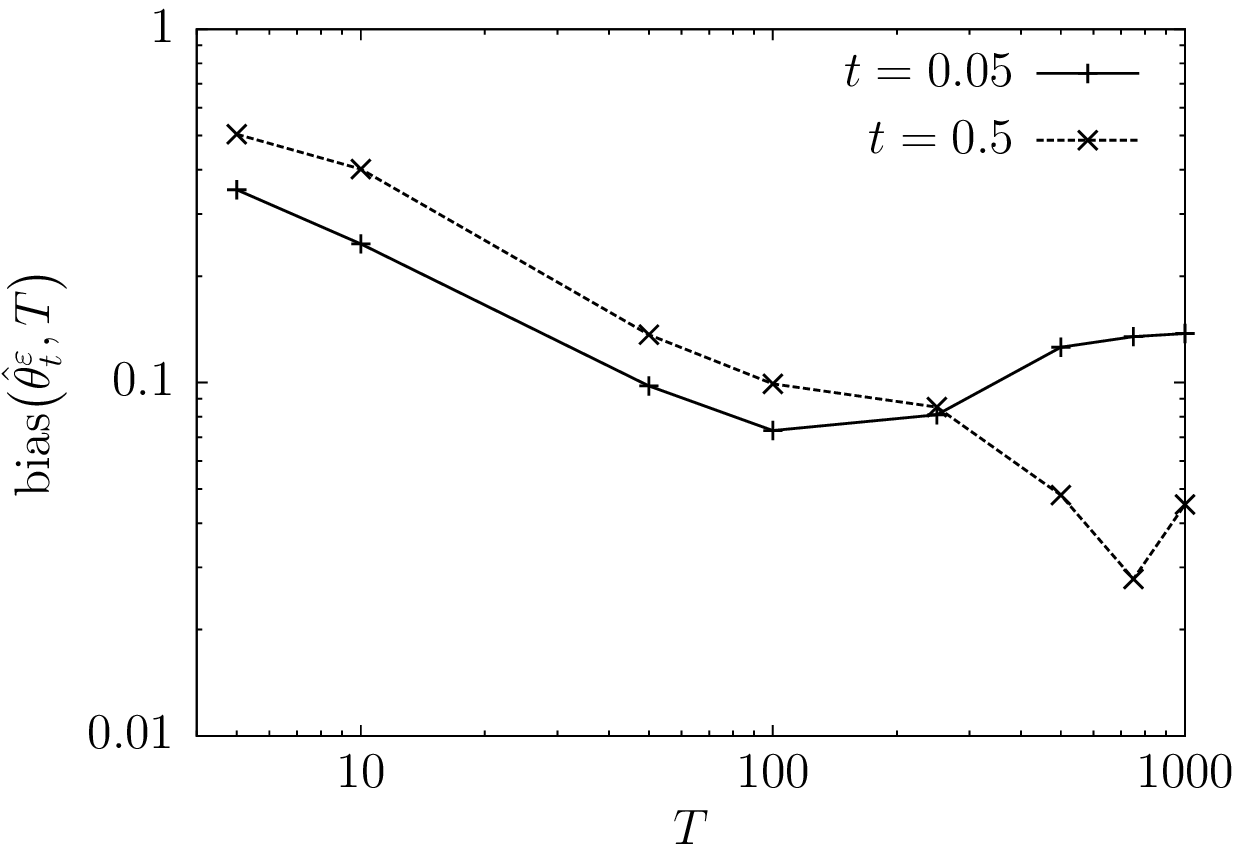}
    \caption{bias}
    \label{figure:langevin:est:bias}
  \end{subfigure}
  \quad
  \begin{subfigure}[b]{0.465\textwidth}
    \centering
    \includegraphics[width=\textwidth]{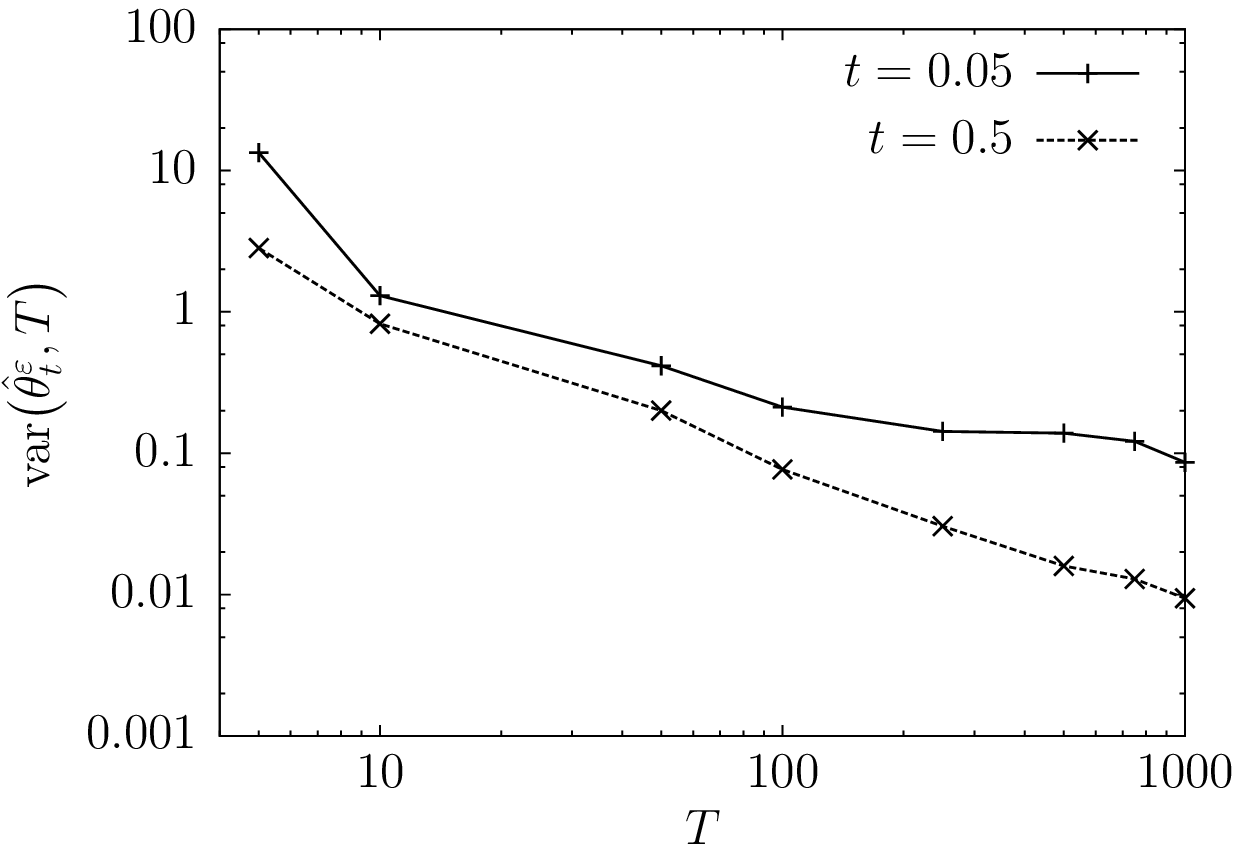}
    \caption{variance}
    \label{figure:langevin:est:var}
  \end{subfigure}
  \caption{Bias and variance of the estimation procedure as functions
    of $T$ for $t\in\{0.05,0.5\}$ and $\varepsilon =
    0.1$. Expectations were approximated as described in the text.}
  \label{figure:langevin:est}
\end{figure}

We proceed by numerically studying the bias and the variance of the
estimation procedure for a fixed time $t$ as functions of the length of the
time interval $T$. To this end we use $M$ independent Brownian motions in
Eq.~\eqref{eq:eps:sde:langevin1d} to generate an ensemble of independent
trajectories $X_1^\varepsilon, X_2^\varepsilon, \dots, X_M^\varepsilon$, each
on the time interval $[0,T]$. Applying the estimation procedure to every such
time series $X_k^\varepsilon$ yields an estimated value, which we denote by
$\hat\theta_t(X_k^\varepsilon,T)$ to emphasize the dependency on the $k$-th
time series and on the final time $T$. Using these estimated values we
approximate expectations by ensemble averages to define the bias and
variance. Specifically, let $E_M^T(\hat\theta_t^\varepsilon) :=
\frac{1}{M}\sum_{k=1}^M\hat\theta_t(X_k^\varepsilon,T)$ be the average of
these estimated values. Then we use
\begin{equation*}
  \bias(\hat\theta_t^\varepsilon,T) :=  {\bigl\Vert E_M^T(\hat\theta_t^\varepsilon) - \theta\bigr\Vert}_2
  \approx {\bigl\Vert \E\bigl(\hat\theta_t(X^\varepsilon,T)\bigr) - \theta\bigr\Vert}_2\;,
\end{equation*}
to quantify the (absolute) bias and as a measure of the variance we
use
\begin{equation*}
  \var(\hat\theta_t^\varepsilon,T) := \frac{1}{M-1}\sum_{i=1}^n\sum_{k=1}^M {\Bigl(e_i\cdot \bigl(\hat\theta_t(X_k^\varepsilon,T) -E_M^T(\hat\theta_t^\varepsilon)\bigr)\Bigr)}^2\approx \sum_{i=1}^n \var\bigl(e_i\cdot \hat\theta_t(X^\varepsilon,T)\bigr)\;,
\end{equation*}
with $e_i$, $1\le i\le n$, denoting the canonical basis vectors of
$\R^n$. In other words $\var(\hat\theta_t^\varepsilon,T)$ is simply an
approximation of the trace of the covariance
matrix. Fig.~\ref{figure:langevin:est} shows the behavior the
estimation procedure's bias and variance as functions of $T$ using
$M=100$ independent Brownian motions for $t\in\{0.05,0.5\}$ and
$\varepsilon = 0.1$.
One observes that the variance $\var(\hat\theta_t^\varepsilon,T)$
(Fig.\ \ref{figure:langevin:est:var}) decreases to zero as $T$
increases, with slightly different rates for the different values of
$t$. Conversely, the bias $\bias(\hat\theta_t^\varepsilon,T)$ (Fig.\
\ref{figure:langevin:est:bias}) starts to decrease by increasing $T$
for both values of $t$, however, after some value of $T$, the bias
approaches a limiting value of approximately $0.15$ for $t=0.05$ and
fluctuates around $0.04$ for $t=0.5$. This fluctuation persists even
for $T>1000$ (not shown here) and are mainly due to the error induced
by approximating an expectation via an ensemble average of size
$M=100$, which becomes visible in this logarithmic scaling. The fact
that the bias approaches a non-zero limiting value is not surprising
(and in fact in agreement with the theoretical results), as one
expects that the estimated value approaches the true value $\theta$,
as $T\rightarrow\infty $, plus an $\mathcal{O}(\varepsilon)$ error due
to the multiscale effects in the data $X^\varepsilon$. Furthermore, we
note that both values of $t$ correspond to estimated values in
Fig.~\ref{figure:langevin} which have a considerable relative error,
where the relative error for $t=0.5$ is significantly smaller that the
one for $t=0.05$, hence explaining the different limiting values in
Fig.\ref{figure:langevin:est:bias} as the constant of the
$\mathcal{O}(\varepsilon)$ error is $t$ dependent.

\subsection{Fast deterministic chaos}
\label{sec:numerics:fastchaos}
We consider an ODE driven by the time rescaled Lorenz equations:
\begin{subequations}
\begin{align}
  \frac{dX^\varepsilon}{dt} & = \alpha\bigl(X^\varepsilon - {(X^\varepsilon )}^3\bigr) + \frac{\lambda}{\varepsilon} Y_2^\varepsilon\;,\label{eq:num:chaos:fastslow:slow}\\
  \frac{dY_1^\varepsilon}{dt} & = \frac{10}{\varepsilon^2}(Y_2^\varepsilon-Y_1^\varepsilon)\;,\label{eq:num:chaos:fastslow:fast1}\\
  \frac{dY_2^\varepsilon}{dt} & = \frac{1}{\varepsilon^2}(28Y_1^\varepsilon - Y_2^\varepsilon - Y_1^\varepsilon Y_3^\varepsilon)\;,\label{eq:num:chaos:fastslow:fast2}\\
  \frac{dY_3^\varepsilon}{dt} & = \frac{1}{\varepsilon^2}\Big(Y_1^\varepsilon Y_2^\varepsilon - \frac{8}{3}Y_3^\varepsilon\Bigr)\;.\label{eq:num:chaos:fastslow:fast3}
\end{align}
\label{eq:num:chaos:fastslow}%
\end{subequations}
Equations of this form have been used as a deterministic climate toy
model, see \cite{Mitchell2012} for instance. Our aim is to obtain a
stochastic coarse-grained model from observations of
\eqref{eq:num:chaos:fastslow}. It is known that, as
$\varepsilon\rightarrow 0$, the slow component $X^\varepsilon$ of
Eq.~\eqref{eq:num:chaos:fastslow} converges weakly in $C([0,T],\R)$ to
the solution of the homogenized equation~\cite{Melbourne2011}
\begin{equation}
  dX = A\bigl(X-X^3\bigr)\,dt + \sqrt{\sigma}\,dW_t\;.\label{eq:num:chaos:effective}
\end{equation}
In Eq.~\eqref{eq:num:chaos:effective} the true parameter values are
$A=\alpha$ and the diffusion coefficient $\sigma$ is given in terms of
the Green--Kubo formula \cite{Pavliotis2008book}
\begin{equation}
  \sigma = 2\lambda^2\int_0^\infty\lim_{T\rightarrow\infty}\frac{1}{T}\int_0^T Y_2^{\varepsilon=1}(s) Y_2^{\varepsilon=1}(s+t)\,ds\,dt\;.
  \label{eq:num:chaos:sigma}
\end{equation}
Obtaining a value for $\sigma$ directly from
Eq.~\eqref{eq:num:chaos:sigma} is computationally challenging so that
the parametric estimation problem of $\sigma$ from observations of
Eq.~\eqref{eq:num:chaos:fastslow:slow} arises naturally for this
model, even without the connection to data-driven coarse-graining
methodologies.

\begin{figure}[t]
  \centering
    \includegraphics[width=0.465\textwidth]{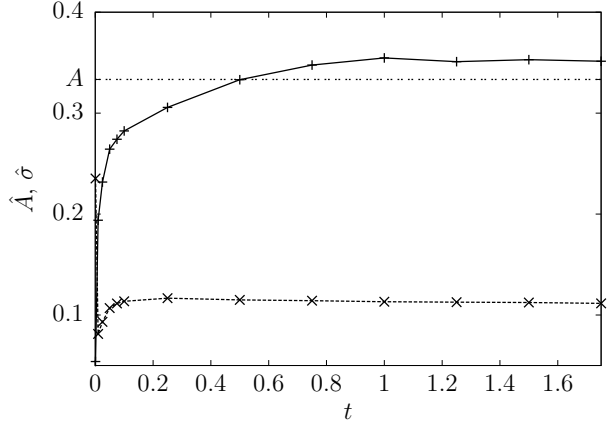}
    \caption[]{Parameter estimates $\hat{A}$~($+$) and
      $\hat{\sigma}$~($\times$) for \eqref{eq:num:chaos:effective}
      based on observations $X^\varepsilon$ of
      \eqref{eq:num:chaos:fastslow} with $\alpha = 1/3$, $\lambda =
      2/45$, and $\varepsilon = 10^{-1}$.}
  \label{figure:fastchaos}
\end{figure}

To estimate both the drift coefficient $A$ and the diffusion
coefficient $\sigma$ (i.e.\ $n=2$) in
Eq.~\eqref{eq:num:chaos:effective}, a self-evident choice for the
functions in Eq.~\eqref{eq:parameterization} is $f_1(x) = x-x^3$,
$f_2(x) = 0 = g_1(x)$, and $g_2(x) = 1$, where the true parameter
vector is $\theta = (A,\sigma)^T$. To generate the time series we
numerically integrate the multiscale system of ODEs
\eqref{eq:num:chaos:fastslow} with $\alpha = 1/3$, $\lambda = 2/45$,
and $\varepsilon=0.1$ on $[0,5000]$ with initial conditions
$X^\varepsilon(0) = 1$, $Y^\varepsilon(0) = (1,1,1)^T$.  For these
parameter choices (mainly the value of $\varepsilon$) the ODE system
\eqref{eq:num:chaos:fastslow} is only marginally stiff and we thus
solve it using a fourth order Runge--Kutta scheme with step size
$h=0.001$. Since there is no exact value for $\sigma$ in
Eq.~\eqref{eq:num:chaos:effective} available, we cannot compute the
relative error of the estimated parameter vector
$\hat{\theta}_t^\varepsilon$. Instead Fig.~\ref{figure:fastchaos}
illustrates both estimated values $\hat{A}$ and $\hat{\sigma}$ as
functions of $t$ directly. One finds that the estimated drift
parameter $\hat{A}$ is strongly biased for very small values of
$t$. Increasing $t$ reduces the bias significantly and the estimated
value approaches the true value (dashed line), only with minor
fluctuations. In fact, the relative error is smaller than $6\%$ for
$t\ge 0.5$. The estimated diffusion coefficient $\hat{\sigma}$ shows
qualitatively the same behavior. Specifically, by increasing $t$ the
estimated value seems to approach a limiting value. In fact, averaging
over the obtained estimated values for $t\ge 0.5$ (i.e.\ the region
for which $\hat{A}$ is accurate), one finds $\hat{\sigma}\approx
0.113$ with minor fluctuations (standard deviation $ \approx
0.002$). This value of $\hat{\sigma}$ is in very good agreement with
those reported in the literature~\cite{Givon2004,Krumscheid2013},
albeit marginally smaller. In fact, the relative error between the
obtained value here and the value reported in \cite{Krumscheid2013} is
around $6\%$.

\subsection{Large Hamiltonian systems: The Kac--Zwanzig model}
\label{sec:numerics:kaczwanzig}
Here we apply our methodology to the variant of the Kac--Zwanzig model
studied in \cite{Kupferman2002}. Specifically, we consider the case
where one distinguished particle, with coordinate $Q_M$ and momentum
$P_M$, moves in an one-dimensional potential $V$ and interacts with
$M\in\N$ heat bath particles. Let the heat bath particles be described
by coordinates $q\equiv (q_1,\dots , q_M)^T\in\R^M$ and momenta
$p\equiv (p_1,\dots , p_M)^T\in\R^M$. Then we consider the Hamiltonian
\begin{equation*}
  H(P_M,Q_M,p,q) := \frac{1}{2}{P_M}^2 + V(Q_M) + \frac{1}{2}\sum_{j=1}^M\frac{{p_j}^2}{m_j}
  + \frac{1}{2}\sum_{j=1}^Mk_j{(q_j-Q_M)}^2\;.
\end{equation*}
That is, the $j$-th heat bath particle with mass $m_j$ acts on the
distinguished particle as a linear spring with stiffness constant $k_j$. The
interaction with the bath is governed by the following ${2(M+1)}$-dimensional
system of ODEs
\begin{subequations}
\begin{align}
  \frac{dQ_M}{dt} &= P_M\;,\qquad \frac{dP_M}{dt} = \sum_{j=1}^Mk_j{(q_j-Q_M)} - V'(Q_M)\;,\\
  \frac{dq_j}{dt} &= \frac{p_j}{m_j}\;,\qquad \frac{dp_j}{dt} = -k_j{(q_j-Q_M)}\;,\quad j=1,\dots, M\;.
\end{align}
\label{eq:kaczwanzig:full:model}%
\end{subequations}
The initial condition for this system are $Q_M(0)=Q_0$, $P_M(0)=P_0$,
$q_j(0)=q_{j,0}$, and $p_j(0)=p_{j,0}$. We assume that the initial
conditions for the heat bath particles are in equilibrium. That is, we
assume that the $2M$-dimensional vector of initial conditions for the
particles in the heat bath (positions and momenta) is randomly
distributed according to a Gibbs distribution with density
proportional to $\exp(-\beta H)$, conditioned on $(Q_0,P_0)$. Here
$\beta >0$ denotes the inverse temperature. Under these conditions it
is possible to derive a coarse-grained model for the distinguished
particle; see e.g.\ \cite{Kupferman2002,Kupferman2004} and the
references therein for details.

The precise form of the coarse-grained model depends mainly on the
chosen values for the spring constants $k_j$ and the particles' mass
$m_j$, $1\le j\le M$. Here we borrow Example $7.3$ from
\cite{Givon2004}. Let $\alpha \in (0,1)$ and define $\omega_j =
M^\alpha\eta_j$, where ${(\eta_j)}_{1\le j\le M}$ is an identically
and independently distributed sequence of random variables with
$\eta_1\sim \mathcal{U}(0,1)$. Moreover, we set $k_j = 2\alpha
M^\alpha/\bigl(\pi({\alpha}^2 + {\omega_j}^2)M \bigr)$ and $m_j =
k_j/{\omega_j}^2$. Then, as $M\rightarrow\infty$, the process $Q_M$
solving the full model \eqref{eq:kaczwanzig:full:model} converges
weakly in $C^2([0,T];\R)$ to the process $Q$ which is the solution of
the stochastic integro-differential equation
\begin{equation*}
  \ddot{Q}(t) + V'\bigl(Q(t)\bigr) + \int_0^t e^{-\alpha\vert t-s\vert}\dot{Q}(s)\,ds = Z(t)\;, 
\end{equation*}
where $Z$ denotes the Ornstein--Uhlenbeck process solving $dZ =
-\alpha Z\,dt + \sqrt{2\alpha/\beta}\,dW$. By introducing an auxiliary
variable, it is possible to convert the integro-differential equation
with nonlocal memory into a Markov process. Specifically, the limiting
process $Q$ is equivalently given as the solution of the augmented
system of SDEs
\begin{subequations}
    \begin{align}
      \frac{dQ}{dt} &= P\;,\label{eq:kaczwanzig:coarse:model:position}\\
      \frac{dP}{dt} &= S - V'(Q)\;,\label{eq:kaczwanzig:coarse:model:momentum}\\
      dS &= (\mu S- P)\,dt + \sqrt{2\sigma}\,dW\;,\label{eq:kaczwanzig:coarse:model:memory}
    \end{align}
\label{eq:kaczwanzig:coarse:model}%
\end{subequations}
where the auxiliary variable $S$ embodies the memory effects due to
the heat bath interactions. Consequently, the coarse-grained model
associated to the ${2(M+1)}$-dimensional Hamiltonian system
\eqref{eq:kaczwanzig:full:model} is given in form of a $3$-dimensional
stochastic system.  The limiting parameters in the coarse-grained
model \eqref{eq:kaczwanzig:coarse:model} are given by $\mu = -\alpha$
and $\sigma = \alpha/\beta$, where we recall that $\beta$ is the
inverse temperature.

The goal now is to estimate $\mu$ and $\sigma$ in
\eqref{eq:kaczwanzig:coarse:model:memory} from observations in form of
a single time series of $(Q_M,P_M)$.  Although the coarse-grained
model \eqref{eq:kaczwanzig:coarse:model} is three-dimensional, we can
use a slightly modified procedure of the one derived in Section
\ref{sec:estimator} for one-dimensional models, since we are concerned
with identifying parameters in only one of the equations in
\eqref{eq:kaczwanzig:coarse:model}, namely in
\eqref{eq:kaczwanzig:coarse:model:memory}. Using It{\^o}'s formula for
\eqref{eq:kaczwanzig:coarse:model} with the function $\phi(s) = s +
s^2$, which only depends on $s$, we find
\begin{equation*}
  \E\Bigl(\phi\bigl(S_\xi(t)\bigr)\Bigr) - \phi(\xi) + \int_0^t\E\Bigl(P_{P_0}(\tau)\phi'\bigl(S_\xi(\tau)\bigr)\Bigr)\,d\tau
  = \int_0^t\E\Bigl((\mathcal{L}_0\phi)\bigl(S_\xi(\tau)\bigr)\Bigr)\,d\tau\;,
\end{equation*}
with $(\mathcal{L}_0\varphi)(s) := \mu s \tfrac{d}{ds}\varphi(s) +
\sigma\tfrac{d^2}{ds^2}\varphi(s)$. This is an estimating equation
like \eqref{eq:ito:est:phi:param} and we thus only have to modify the
definition of the term $b_c$ in \eqref{eq:fun:form:param:single} to
account for the dependency of
\eqref{eq:kaczwanzig:coarse:model:memory} on the process $P$ (the
integral term on the left-hand side above). The rest of the procedure
follows as in Section \ref{sec:estimator}.  In fact, we select the
functions in parametrization \eqref{eq:parameterization} with $n=2$ as
$f_1(x) = x$, $f_2(x) = 0 = g_1(x)$, and $g_2(x) = 2$, where the true
parameter vector is $\theta = (\mu,\sigma)^T$.

It is important to stress that we wish to estimate $\mu$ and $\sigma$
in \eqref{eq:kaczwanzig:coarse:model:memory}, but that we do not
observe the process $S$ directly: unlike for $Q$ and $P$ where we
observe $Q_M$ and $P_M$ which converge to $Q$ and $P$, respectively,
we do not have access to such a process for $S$. We only have
observations of $(Q_M,P_M)$ from the full Hamiltonian system
\eqref{eq:kaczwanzig:full:model} with sampling time $h$.  In the
absence of model misspecification (i.e., when observing $Q,P$ directly
and not just $Q_M,P_M$ instead), this problem is typically associated
to hidden Markov model techniques as $S$ is unobserved (i.e.\ hidden);
see e.g. \cite{Cappe2005}.  Here we consider a simple approximation to
reconstruct the unobserved process $S$, which we will need in the
estimation procedure. To this end we use the observations we have in
\eqref{eq:kaczwanzig:coarse:model:momentum} with a first order finite
difference approximation:
\begin{equation*}
S_M(t) := \frac{P_M(t+h)-P_M(t)}{h} + V'\bigl(Q_M(t)\bigr)\;.
\end{equation*}
We remark that, in principle, we can apply our methodology even if
only $Q_M$ is observed but not $P_M$. In that case one has to use both
Eq.~\eqref{eq:kaczwanzig:coarse:model:position} and
Eq.~\eqref{eq:kaczwanzig:coarse:model:momentum} with finite difference
approximations to obtain suitable approximations of $P$ and $S$. These
finite difference approximation ideas have also been used in
\cite{Kupferman2004} within a customized maximum likelihood framework
for the Kac--Zwanzig model. However, in their study the authors had to
chose $h$ sufficiently large as otherwise the parameter estimation
performed poorly due to the presence of multiscale effects. Here we
are not restricted by the multiscale structure of the problem.

\begin{figure}[t]
  \centering
    \includegraphics[width=0.465\textwidth]{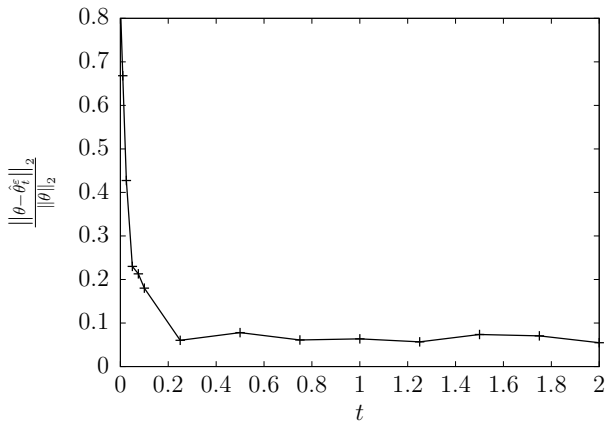}
    \caption[]{Relative error of the estimated parameter vector
      $\hat{\theta}_{t}^\varepsilon$ for
      \eqref{eq:kaczwanzig:coarse:model:memory} based on observations
      $(Q_M,P_M)$ of \eqref{eq:kaczwanzig:full:model} with $\alpha
      =1/2 2$, $\beta=1$, and $M\equiv \varepsilon^{-1} = 5000$.}
  \label{figure:hamiltonian}
\end{figure}

For the numerical example we consider the case where the distinguished
particle moves in a quartic potential, i.e.\ $V(x) = -x^2/2 + x^4/4$.
Moreover, we use $M = 5000$ heat bath particles and set $\alpha = 1/2$
and $\beta= 1$. To obtain a time series for $(Q_M,P_M)$ of the full
Hamiltonian system \eqref{eq:kaczwanzig:full:model} on $[0,1000]$, we
approximate it via a semi-implicit Euler scheme with time step
$h=10^{-3}$ started at $Q_0 = 1$, $P_0 =
0$. Fig.~\ref{figure:hamiltonian} depicts the relative error
of the estimated parameter vector $\hat{\theta}_t^\varepsilon$ as a
function of $t$, with the understanding that $\varepsilon \equiv
M^{-1}$. Similar to the previous examples, one also observes here that
it is possible to obtain accurate estimates once the value of $t$ is
sufficiently large. The relative error fluctuates closely around $6\%$
for $t\ge 0.2$ and can be reduced even further by increasing $M$ (not
shown). To asses the obtained accuracy of $6\%$, we mention that in
\cite{Kupferman2004}, as noted above, a specifically customized
maximum likelihood method has been used to estimate parameters in a
related problem and the authors report relative errors of $5$--$15\%$.

\subsection{Deterministic Brownian motion}
\label{sec:numerics:detbm}
In \cite{Mackey2006} an Ornstein--Uhlenbeck process is constructed
within a completely deterministic framework as an appropriate limit
process of a chaotic dynamical system. Specifically, consider the
position $X^\varepsilon$ and the velocity $V^\varepsilon$ of the
dynamical system
\begin{subequations}
\begin{align}
  \frac{dX^\varepsilon}{dt} &= V^\varepsilon\;,\label{eq:detBM:full:dx}\\
  \frac{dV^\varepsilon}{dt} &= -\gamma V^\varepsilon + \eta_\varepsilon\;,\label{eq:detBM:full:dv}
\end{align}
\label{eq:detBM:full}%
\end{subequations}
with a deterministic perturbation
$\eta_\varepsilon\equiv\eta_\varepsilon(t)$ in the velocity
variable. Here we use
\begin{equation}
  \eta_\varepsilon(t) = \sqrt{\varepsilon}\sum_{l=0}^{\infty}\zeta(t_l)\delta(t-t_l)\;,\label{eq:detBM:perturb}
\end{equation}
so that the derivative of the velocity variable $V^\varepsilon$
experiences small ``kicks'' at times $t_0,t_1,\dots$, where $t_l =
l\varepsilon$. Here $\zeta$ is a discrete time dynamical system of the
form $\zeta(t_{l+1}) = \Phi(\zeta(t_{l}))$ and the function $\Phi$ is
chosen such that the dynamical system $\zeta$ exhibits a strongly
chaotic behavior. In our numerical example we used $\Phi(y) :=
\cos{(3\arccos(y))}$. Based on this dynamical system $\zeta$ the
perturbation $\eta_\varepsilon$ in Eq.~\eqref{eq:detBM:perturb} is
fixed and it follows from the results in \cite{Mackey2006} that the
solution $(X^\varepsilon,V^\varepsilon)$ of the chaotic deterministic
system \eqref{eq:detBM:full} converges weakly in $C([0,T],\R)$, as
$\varepsilon\rightarrow 0$, to an Ornstein--Uhlenbeck process $(X,V)$,
which solves
\begin{subequations}
\begin{align}
  dX &= V\,dt\;,\label{eq:detBM:coarse:dx}\\
  dV &= -\gamma V\,dt + \sqrt{\sigma}\,dW_t\;,\label{eq:detBM:coarse:dv}
\end{align}
\label{eq:detBM:coarse}%
\end{subequations}
where the diffusion coefficient is $\sigma = 1/2$.

\begin{figure}[t]
  \centering
    \includegraphics[width=0.465\textwidth]{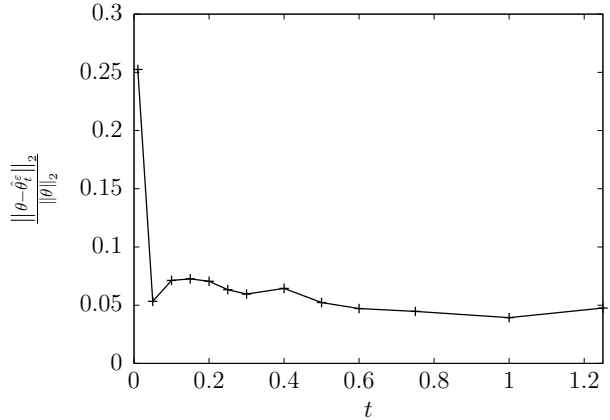}
    \caption[]{Relative error of the estimated parameter vector
      $\hat{\theta}_{t}^\varepsilon$ for \eqref{eq:detBM:coarse:dv}
      based on observations $X^\varepsilon$ of \eqref{eq:detBM:full}
      with $\gamma = 1$ and $\varepsilon = 0.1$.}
  \label{figure:detbm}
\end{figure}

We now aim for estimating both $\gamma$ and $\sigma$ in
\eqref{eq:detBM:coarse:dv} based only on one long trajectory of
observations of the position variable $X^\varepsilon$ solving
\eqref{eq:detBM:full:dx}. That is, we do not observe $V^\varepsilon$
solving \eqref{eq:detBM:full:dv} directly.  Instead we will, as in
Section \ref{sec:numerics:kaczwanzig}, compute an approximation
$\tilde{V}^\varepsilon$ based on a finite difference approximation in
\eqref{eq:detBM:coarse:dx} first, i.e.\ we set
$\tilde{V}^\varepsilon(t_l):=(X^\varepsilon(t_{l+1})-X^\varepsilon(t_{l}))/\varepsilon$,
recalling that $t_{l+1}-t_{l} = \varepsilon$. Based on this
approximate trajectory we can then directly apply the procedure
introduced in Section \ref{sec:estimator} to estimate both $\gamma$
and $\sigma$ in \eqref{eq:detBM:coarse:dv}, since the velocity SDE is
independent of the position. Therefore ($n=2$) and we select the
functions $f_1(x) = x$, $f_2(x) = 0 = g_1(x)$, and $g_2(x) = 1$ in
\eqref{eq:parameterization}, corresponding to the true parameter
vector $\theta = (-\gamma,\sigma)^T$. A time series of $X^\varepsilon$
on $[0,1000]$ with sampling rate $h=0.01$ is obtained by solving the
perturbed system \eqref{eq:detBM:full} with $X^\varepsilon(0) =
-0.15$, $V^\varepsilon(0)=-0.53$, $\gamma = 1$ and $\varepsilon =0.1$.
Fig.~\ref{figure:detbm} shows
the relative error of the estimated parameter vector
$\hat{\theta}_t^\varepsilon$ as a function of $t$. Increasing $t$ yields very
accurate estimates with the relative error fluctuating around $5\%$ for $t\ge
0.5$.

%
%
\section{Conclusion}
\label{sec:conclusion}
In this paper we have introduced a novel numerical/statistical
procedure which allows us to estimate parameters in coarse-grained
models based on partial observations of a corresponding multiscale
system. For such systems commonly used estimators, such as the maximum
likelihood estimator, are known to be biased. Our approach is based on
our previous study in~\cite{Krumscheid2013} where it was assumed that
an ensemble of short trajectories for multiple initial conditions is
available.  Here we generalize and appropriately extend the work
presented in \cite{Krumscheid2013} to the practically relevant setting
where only one (long) time series is available. In fact, the examples
presented demonstrate that the developed inference method yields
accurate approximations of the parameters in coarse-grained models
based on a time series of the ``slow'' component of a multiscale
system. The examples range from coarse-grained models where the
associated multiscale system is stochastic to coarse-grained models
for fully deterministic multiscale systems. We believe that this
selection of examples highlights the broad occurrence of data-driven
coarse-graining problems and thus the necessity for appropriate
inference techniques which are robust against multiscale effects in
the observation, as the one introduced here.

The illustrated robustness of an inference technique against
multiscale effects (i.e.\ against perturbations that are small the
weak sense) appears to be novel but also has important implications in
practice. Specifically, it significantly widens the range of
applications to problems where the observed process is not necessarily
a diffusion process. For instance, it covers the case of non-Markovian
processes that can be approximated by a Markov process in an augmented
state space (see, e.g., \cite[Ch.~$8.2$]{Pavliotis2014_book}) and the
concept of diffusion approximation (see, e.g.,
\cite[Ch.~$7$]{Ethier1986}), both allowing for non-negligible
deviations of the observation process from the assumed SDE model.

The focus of our study was on demonstrating that the introduced
methodology can accurately infer parameters in coarse-grained models
from a time series, either stochastic or chaotic, of a multiscale
system. Clearly there are still many challenges that remain to be
addressed. One of them is the rigorous analysis of the algorithm to
understand its asymptotic properties, but also to explore its
limitations. Some first results concerning the convergence properties
of the methodology have already been obtained and will be presented in
\cite{Krumscheid_pre}. A closely related and important avenue of
future efforts is the study of the asymptotic distribution of the
estimators, which in turn can be used to guide the construction of
asymptotic confidence intervals for the estimated values.  Another
interesting topic, which is relevant in many applications, is the
study of additional observation error. That is, one only observes a
contaminated version $\tilde{X}^\varepsilon$ of $X^\varepsilon$:
\begin{equation*}
  \tilde{X}^\varepsilon(t_k) = X^\varepsilon(t_k) +  \eta(t_k)\;,
\end{equation*}
for any $k\ge 0$, where $\eta$ denotes the observation error. To allow
for this additional contamination, our methodology would have to be
combined with appropriate filtering techniques, a very appealing
prospect. The first results on combining filtering ideas with
parametric inference techniques have been studied in \cite{Cotter2009}
for a particular multiscale problem. See also \cite{Papanicolaou2014}
for more recent work on combining the MLE with filtering techniques
for a class of multiscale problems.  But also conceptually different
approaches to the problem of data-driven coarse-graining appear
worthwhile investigating.  As, for example, most work on data-driven
coarse-graining is based on a frequentist inference approach,
investigating similar questions in a Bayesian approach poses a natural
and interesting perspective as well; see \cite{Nolen2012} for related
work in the context of inverse problems with a multiscale
structure. We shall examine these and related questions in future
studies.

%
%
\section*{Acknowledgements}
We are grateful to the anonymous referees for their insightful
comments and suggestions. We acknowledge financial support from the
Engineering and Physical Sciences Research Council of the UK through
Grants No.\ EP/H034587, EP/J009636, and EP/L020564 and from the
European Research Council via Advanced Grant No.\ 247031.

\bibliography{references}

\begin{thebibliography}{10}
\expandafter\ifx\csname url\endcsname\relax
  \def\url#1{\texttt{#1}}\fi
\expandafter\ifx\csname urlprefix\endcsname\relax\def\urlprefix{URL }\fi
\expandafter\ifx\csname href\endcsname\relax
  \def\href#1#2{#2} \def\path#1{#1}\fi

\bibitem{Chauviere2010}
A.~Chauvi{\`e}re, L.~Preziosi, C.~Verdier (Eds.), Cell Mechanics: {F}rom Single
  Scale-Based Models to Multiscale Modeling, Mathematical \& Computational
  Biology Series, Chapman \& Hall/CRC, 2010.

\bibitem{Majda2008}
A.~J. Majda, C.~Franzke, B.~Khouider, An applied mathematics perspective on
  stochastic modelling for climate, Philos. Trans. R. Soc. Lond. Ser. A Math.
  Phys. Eng. Sci. 366~(1875) (2008) 2429--2455.
\newblock \href {http://dx.doi.org/10.1098/rsta.2008.0012}
  {\path{doi:10.1098/rsta.2008.0012}}.

\bibitem{Culina2010}
J.~Culina, S.~Kravtsov, A.~H. Monahan, Stochastic parameterization schemes for
  use in realistic climate models, J. Atmospheric Sci. 68~(2) (2010) 284--299.
\newblock \href {http://dx.doi.org/10.1175/2010JAS3509.1}
  {\path{doi:10.1175/2010JAS3509.1}}.

\bibitem{Griebel2007}
M.~Griebel, S.~Knapek, G.~W. Zumbusch, Numerical Simulation in Molecular
  Dynamics: Numerics, Algorithms, Parallelization, Applications, Texts in
  Computational Science and Engineering, Springer, 2007.

\bibitem{Fish2009}
J.~Fish, Multiscale Methods: {B}ridging the Scales in Science and Engineering,
  Oxford University Press, 2009.

\bibitem{Huerre1998}
P.~Huerre, M.Rossi, Hydrodynamic instabilities in open flows, in:
  C.~Godr\`eche, P.~Manneville (Eds.), {H}ydrodynamic and {N}onlinear
  {I}nstabilities, Cambridge University Press, 1998, pp. 81--294.

\bibitem{Horstemeyer2010}
M.~F. Horstemeyer, Multiscale modeling: {A} review, in: J.~Leszczynski, M.~K.
  Shukla (Eds.), Practical Aspects of Computational Chemistry, Springer, 2010,
  pp. 87--135.
\newblock \href {http://dx.doi.org/10.1007/978-90-481-2687-3}
  {\path{doi:10.1007/978-90-481-2687-3}}.

\bibitem{Savva2010}
N.~Savva, S.~Kalliadasis, G.~A. Pavliotis, Two-dimensional droplet spreading
  over random topographical substrates, Phys. Rev. Lett. 104~(8) (2010) 084501.
\newblock \href {http://dx.doi.org/10.1103/PhysRevLett.104.084501}
  {\path{doi:10.1103/PhysRevLett.104.084501}}.

\bibitem{Chorin2000}
A.~J. Chorin, O.~H. Hald, R.~Kupferman, Optimal prediction and the
  {M}ori-{Z}wanzig representation of irreversible processes, Proc. Natl. Acad.
  Sci. USA 97~(7) (2000) 2968--2973.
\newblock \href {http://dx.doi.org/10.1073/pnas.97.7.2968}
  {\path{doi:10.1073/pnas.97.7.2968}}.

\bibitem{Turkington2013}
B.~Turkington, An optimization principle for deriving nonequilibrium
  statistical models of {H}amiltonian dynamics, J. Stat. Phys. 152~(3) (2013)
  569--597.
\newblock \href {http://dx.doi.org/10.1007/s10955-013-0778-9}
  {\path{doi:10.1007/s10955-013-0778-9}}.

\bibitem{Venturi2014}
D.~Venturi, G.~E. Karniadakis, Convolutionless {N}akajima--{Z}wanzig equations
  for stochastic analysis in nonlinear dynamical systems, Proc. R. Soc. Lond.
  Ser. A Math. Phys. Eng. Sci. 470~(2166) (2014) 20130754, 20.
\newblock \href {http://dx.doi.org/10.1098/rspa.2013.0754}
  {\path{doi:10.1098/rspa.2013.0754}}.

\bibitem{Pavliotis2007}
G.~A. Pavliotis, A.~M. Stuart, Parameter estimation for multiscale diffusions,
  J. Stat. Phys. 127~(4) (2007) 741--781.
\newblock \href {http://dx.doi.org/10.1007/s10955-007-9300-6}
  {\path{doi:10.1007/s10955-007-9300-6}}.

\bibitem{Krumscheid2013}
S.~Krumscheid, G.~A. Pavliotis, S.~Kalliadasis, Semiparametric drift and
  diffusion estimation for multiscale diffusions, Multiscale Model. Simul.
  11~(2) (2013) 442--473.
\newblock \href {http://dx.doi.org/10.1137/110854485}
  {\path{doi:10.1137/110854485}}.

\bibitem{Pradas2011}
M.~Pradas, D.~Tseluiko, S.~Kalliadasis, D.~T. Papageorgiou, G.~A. Pavliotis,
  Noise induced state transitions, intermittency, and universality in the noisy
  {K}uramoto-{S}ivashinksy equation, Phys. Rev. Lett. 106~(6) (2011) 060602.
\newblock \href {http://dx.doi.org/10.1103/PhysRevLett.106.060602}
  {\path{doi:10.1103/PhysRevLett.106.060602}}.

\bibitem{Pradas2012}
M.~Pradas, G.~A. Pavliotis, S.~Kalliadasis, D.~T. Papageorgiou, D.~Tseluiko,
  Additive noise effects in active nonlinear spatially extended systems,
  European J. Appl. Math. 23~(5) (2012) 563--591.
\newblock \href {http://dx.doi.org/10.1017/S0956792512000125}
  {\path{doi:10.1017/S0956792512000125}}.

\bibitem{Schmuck2013}
M.~Schmuck, M.~Pradas, S.~Kalliadasis, G.~A. Pavliotis, New stochastic mode
  reduction strategy for dissipative systems, Phys. Rev. Lett. 110~(24) (2013)
  244101.
\newblock \href {http://dx.doi.org/10.1103/PhysRevLett.110.244101}
  {\path{doi:10.1103/PhysRevLett.110.244101}}.

\bibitem{Pavliotis2008book}
G.~A. Pavliotis, A.~M. Stuart, Multiscale Methods: Averaging and
  Homogenization, Springer, 2008.

\bibitem{PrakasaRao1999}
B.~L.~S. {Prakasa Rao}, Statistical Inference for Diffusion Type Processes,
  Vol.~8 of Kendall's Library of Statistics, Arnold, 1999.

\bibitem{Kutoyants2004}
Y.~A. Kutoyants, Statistical Inference for Ergodic Diffusion Processes,
  Springer, 2004.
\newblock \href {http://dx.doi.org/10.1007/978-1-4471-3866-2}
  {\path{doi:10.1007/978-1-4471-3866-2}}.

\bibitem{Liptser2010}
R.~S. Liptser, A.~N. Shiryaev, Statistics of Random Processes: {I}. General
  Theory, 2nd Edition, Stochastic Modelling and Applied Probability Series,
  Springer, 2010, translated by A. B. Aries.

\bibitem{Papavasiliou2009}
A.~Papavasiliou, G.~A. Pavliotis, A.~M. Stuart, Maximum likelihood drift
  estimation for multiscale diffusions, Stochastic Process. Appl. 119 (2009)
  3173--3210.
\newblock \href {http://dx.doi.org/10.1016/j.spa.2009.05.003}
  {\path{doi:10.1016/j.spa.2009.05.003}}.

\bibitem{Zhang2005}
L.~Zhang, P.~A. Mykland, Y.~A{\"{\i}}t-Sahalia, A tale of two time scales:
  determining integrated volatility with noisy high-frequency data, J. Amer.
  Statist. Assoc. 100~(472) (2005) 1394--1411.
\newblock \href {http://dx.doi.org/10.1198/016214505000000169}
  {\path{doi:10.1198/016214505000000169}}.

\bibitem{Azencott2010}
R.~Azencott, A.~Beri, I.~Timofeyev, Adaptive sub-sampling for parametric
  estimation of gaussian diffusions, J. Stat. Phys. 139~(6) (2010) 1066--1089.
\newblock \href {http://dx.doi.org/10.1007/s10955-010-9975-y}
  {\path{doi:10.1007/s10955-010-9975-y}}.

\bibitem{Azencott2011}
R.~Azencott, A.~Beri, I.~Timofeyev, Parametric estimation of stationary
  stochastic processes under indirect observability, J. Stat. Phys. 144~(1)
  (2011) 150--170.
\newblock \href {http://dx.doi.org/10.1007/s10955-011-0253-4}
  {\path{doi:10.1007/s10955-011-0253-4}}.

\bibitem{Cotter2009}
C.~J. Cotter, G.~A. Pavliotis,
  \href{http://projecteuclid.org/euclid.cms/1264434134}{Estimating eddy
  diffusivities from noisy {L}agrangian observations}, Commun. Math. Sci. 7~(4)
  (2009) 805--838.
\newline\urlprefix\url{http://projecteuclid.org/euclid.cms/1264434134}

\bibitem{Olhede2009}
S.~C. Olhede, A.~M. Sykulski, G.~A. Pavliotis, Frequency domain estimation of
  integrated volatility for {I}t\^o processes in the presence of
  market-microstructure noise, Multiscale Model. Simul. 8~(2) (2009) 393--427.
\newblock \href {http://dx.doi.org/10.1137/090756363}
  {\path{doi:10.1137/090756363}}.

\bibitem{Crommelin2011}
D.~T. Crommelin, E.~Vanden-Eijnden, Diffusion estimation from multiscale data
  by operator eigenpairs, Multiscale Model. Simul. 9~(4) (2011) 1588--1623.
\newblock \href {http://dx.doi.org/10.1137/100795917}
  {\path{doi:10.1137/100795917}}.

\bibitem{Crommelin2012}
D.~Crommelin, Estimation of space-dependent diffusions and potential landscapes
  from non-equilibrium data, J. Stat. Phys. 149~(2) (2012) 220--233.
\newblock \href {http://dx.doi.org/10.1007/s10955-012-0597-4}
  {\path{doi:10.1007/s10955-012-0597-4}}.

\bibitem{Spiliopoulos2013}
K.~Spiliopoulos, A.~Chronopoulou, Maximum likelihood estimation for small noise
  multiscale diffusions, Stat. Inference Stoch. Process. 16~(3) (2013)
  237--266.
\newblock \href {http://dx.doi.org/10.1007/s11203-013-9088-8}
  {\path{doi:10.1007/s11203-013-9088-8}}.

\bibitem{Imkeller2013}
P.~Imkeller, N.~{Sri Namachchivaya}, N.~Perkowski, H.~C. Yeong, Dimensional
  reduction in nonlinear filtering: {A} homogenization approach, Ann. Appl.
  Probab. 23~(6) (2013) 2290--2326.
\newblock \href {http://dx.doi.org/10.1214/12-AAP901}
  {\path{doi:10.1214/12-AAP901}}.

\bibitem{Zhang2014}
W.~Zhang, J.~C. Latorre, G.~A. Pavliotis, C.~Hartmann, Optimal control of
  multiscale systems using reduced-order models, J. Computational Dynamics
  1~(2) (2014) 279--306.
\newblock \href {http://dx.doi.org/10.3934/jcd.2014.1.279}
  {\path{doi:10.3934/jcd.2014.1.279}}.

\bibitem{Vanden-Eijnden2003}
E.~Vanden-Eijnden,
  \href{http://projecteuclid.org/euclid.cms/1118152078}{Numerical techniques
  for multi-scale dynamical systems with stochastic effects}, Commun. Math.
  Sci. 1~(2) (2003) 385--391.
\newline\urlprefix\url{http://projecteuclid.org/euclid.cms/1118152078}

\bibitem{E2005}
W.~E, D.~Liu, E.~Vanden-Eijnden, Analysis of multiscale methods for stochastic
  differential equations, Comm. Pure Appl. Math. 58~(11) (2005) 1544--1585.
\newblock \href {http://dx.doi.org/http://dx.doi.org/10.1002/cpa.20088}
  {\path{doi:http://dx.doi.org/10.1002/cpa.20088}}.

\bibitem{Theodoropoulos2000}
C.~Theodoropoulos, Y.-H. Qian, I.~G. Kevrekidis,
  \href{http://www.jstor.org/stable/123274}{"coarse" stability and bifurcation
  analysis using time-steppers: {A} reaction-diffusion example}, Proc. Natl.
  Acad. Sci. USA 97~(18) (2000) 9840--9843.
\newline\urlprefix\url{http://www.jstor.org/stable/123274}

\bibitem{Kevrekidis2003}
I.~G. Kevrekidis, C.~W. Gear, J.~M. Hyman, P.~G. Kevrekidis, O.~Runborg,
  C.~Theodoropoulos,
  \href{http://projecteuclid.org/getRecord?id=euclid.cms/1119655353}{Equation-free,
  coarse-grained multiscale computation: enabling microscopic simulators to
  perform system-level analysis}, Commun. Math. Sci. 1~(4) (2003) 715--762.
\newline\urlprefix\url{http://projecteuclid.org/getRecord?id=euclid.cms/1119655353}

\bibitem{Kevrekidis2004}
I.~G. Kevrekidis, C.~W. Gear, G.~Hummer, Equation-free: {T}he computer-aided
  analysis of complex multiscale systems, AIChE J. 50~(7) (2004) 1346--1355.
\newblock \href {http://dx.doi.org/10.1002/aic.10106}
  {\path{doi:10.1002/aic.10106}}.

\bibitem{Kevrekidis2009}
I.~G. Kevrekidis, G.~Samaey, Equation-free multiscale computation: {A}lgorithms
  and applications, Annu. Rev. Phys. Chem. 60~(1) (2009) 321--344.
\newblock \href {http://dx.doi.org/10.1146/annurev.physchem.59.032607.093610}
  {\path{doi:10.1146/annurev.physchem.59.032607.093610}}.

\bibitem{Karatzas1991}
I.~Karatzas, S.~E. Shreve, Brownian Motion and Stochastic Calculus, 2nd
  Edition, Springer, 1991.
\newblock \href {http://dx.doi.org/10.1007/978-1-4612-0949-2}
  {\path{doi:10.1007/978-1-4612-0949-2}}.

\bibitem{Oksendal2003}
B.~K. {\O}ksendal, Stochastic Differential Equations: An Introduction with
  Applications, Springer, 2003.
\newblock \href {http://dx.doi.org/10.1007/978-3-642-14394-6}
  {\path{doi:10.1007/978-3-642-14394-6}}.

\bibitem{Krumscheid_pre}
S.~Krumscheid, Perturbation-based inference for diffusion processes:
  {O}btaining coarse-grained models from multiscale data, submitted (2014).

\bibitem{Bosq1998}
D.~Bosq, Nonparametric Statistics for Stochastic Processes: Estimation and
  Prediction, 2nd Edition, Vol. 110 of Lecture Notes in Statistics, Springer,
  1998.
\newblock \href {http://dx.doi.org/10.1007/978-1-4612-1718-3}
  {\path{doi:10.1007/978-1-4612-1718-3}}.

\bibitem{Fan2003}
J.~Fan, Q.~Yao, Nonlinear Time Series: Nonparametric and Parametric Methods,
  Springer Series in Statistics, Springer, 2003.
\newblock \href {http://dx.doi.org/10.1007/b97702} {\path{doi:10.1007/b97702}}.

\bibitem{Doukhan1994}
P.~Doukhan, Mixing: Properties and Examples, Springer, 1994.
\newblock \href {http://dx.doi.org/10.1007/978-1-4612-2642-0}
  {\path{doi:10.1007/978-1-4612-2642-0}}.

\bibitem{ChenX2010}
X.~Chen, L.~P. Hansen, M.~Carrasco, Nonlinearity and temporal dependence, J.
  Econometrics 155~(2) (2010) 155--169.
\newblock \href {http://dx.doi.org/10.1016/j.jeconom.2009.10.001}
  {\path{doi:10.1016/j.jeconom.2009.10.001}}.

\bibitem{Nadaraya1964}
E.~A. Nadaraya, On estimating regression, Theory Probab. Appl. 9~(1) (1964)
  141--142.
\newblock \href {http://dx.doi.org/10.1137/1109020}
  {\path{doi:10.1137/1109020}}.

\bibitem{Watson1964}
G.~S. Watson, \href{http://www.jstor.org/stable/25049340}{Smooth regression
  analysis}, Sankhy\=a Ser. A 26~(4) (1964) 359--372.
\newline\urlprefix\url{http://www.jstor.org/stable/25049340}

\bibitem{Cruz-Uribe2002}
D.~Cruz-Uribe, C.~J. Neugebauer,
  \href{http://www.emis.de/journals/JIPAM/article201.html?sid=201}{Sharp error
  bounds for the trapezoidal rule and {S}impson's rule}, JIPAM. J. Inequal.
  Pure Appl. Math. 3~(4) (2002) Article 49, 22.
\newline\urlprefix\url{http://www.emis.de/journals/JIPAM/article201.html?sid=201}

\bibitem{Golub1996}
G.~H. Golub, C.~F. van Loan, Matrix Computations, 3rd Edition, Johns Hopkins
  Studies in the Mathematical Sciences, Johns Hopkins University Press, 1996.

\bibitem{Kroese2011}
D.~P. Kroese, T.~Taimre, Z.~I. Botev, Handbook of Monte Carlo Methods, Wiley
  Series in Probability and Statistics, John Wiley and Sons, 2011.

\bibitem{Mitchell2012}
L.~Mitchell, G.~A. Gottwald, Data assimilation in slow-fast systems using
  homogenized climate models, J. Atmospheric Sci. 69~(4) (2012) 1359--1377.
\newblock \href {http://dx.doi.org/10.1175/JAS-D-11-0145.1}
  {\path{doi:10.1175/JAS-D-11-0145.1}}.

\bibitem{Melbourne2011}
I.~Melbourne, A.M.Stuart, A note on diffusion limits of chaotic skew-product
  flows, Nonlinearity 24~(4) (2011) 1361--1367.
\newblock \href {http://dx.doi.org/10.1088/0951-7715/24/4/018}
  {\path{doi:10.1088/0951-7715/24/4/018}}.

\bibitem{Givon2004}
D.~Givon, R.~Kupferman, A.~M. Stuart, Extracting macroscopic dynamics: model
  problems and algorithms, Nonlinearity 17~(6) (2004) 55--127.
\newblock \href {http://dx.doi.org/10.1088/0951-7715/17/6/R01}
  {\path{doi:10.1088/0951-7715/17/6/R01}}.

\bibitem{Kupferman2002}
R.~Kupferman, A.~M. Stuart, J.~R. Terry, P.~F. Tupper, Long-term behaviour of
  large mechanical systems with random initial data, Stoch. Dyn. 2~(4) (2002)
  533--562.
\newblock \href {http://dx.doi.org/10.1142/S0219493702000571}
  {\path{doi:10.1142/S0219493702000571}}.

\bibitem{Kupferman2004}
R.~Kupferman, A.~M. Stuart, Fitting {SDE} models to nonlinear {K}ac-{Z}wanzig
  heat bath models, Phys. D 199~(3-4) (2004) 279--316.
\newblock \href {http://dx.doi.org/10.1016/j.physd.2004.04.011}
  {\path{doi:10.1016/j.physd.2004.04.011}}.

\bibitem{Cappe2005}
O.~Capp{\'e}, E.~Moulines, T.~Ryd{\'e}n, Inference in Hidden {M}arkov Models,
  Springer Series in Statistics, Springer, 2005.

\bibitem{Mackey2006}
M.~C. Mackey, M.~Tyran-Kami{\'n}ska, Deterministic {B}rownian motion: the
  effects of perturbing a dynamical system by a chaotic semi-dynamical system,
  Phys. Rep. 422~(5) (2006) 167--222.
\newblock \href {http://dx.doi.org/10.1016/j.physrep.2005.09.002}
  {\path{doi:10.1016/j.physrep.2005.09.002}}.

\bibitem{Pavliotis2014_book}
G.~A. Pavliotis, Stochastic Processes and Applications: Diffusion Processes,
  the Fokker-Planck and Langevin Equations, Springer, 2014.

\bibitem{Ethier1986}
S.~N. Ethier, T.~G. Kurtz, Markov Processes: Characterization and Convergence,
  Wiley Series in Probability and Mathematical Statistics, John Wiley \& Sons,
  Inc., 1986.
\newblock \href {http://dx.doi.org/10.1002/9780470316658}
  {\path{doi:10.1002/9780470316658}}.

\bibitem{Papanicolaou2014}
A.~Papanicolaou, K.~Spiliopoulos, Filtering the maximum likelihood for
  multiscale problems, Multiscale Model. Simul. 12~(3) (2014) 1193--1229.
\newblock \href {http://dx.doi.org/10.1137/140952648}
  {\path{doi:10.1137/140952648}}.

\bibitem{Nolen2012}
J.~Nolen, G.~A. Pavliotis, A.~M. Stuart, Multiscale modelling and inverse
  problems, in: O.~L. I.G.~Graham, T.Y.~Hou, R.~Scheichl (Eds.), Numerical
  Analysis of Multiscale Problems, Vol.~83, Springer, 2012, pp. 1--34.
\newblock \href {http://dx.doi.org/10.1007/978-3-642-22061-6}
  {\path{doi:10.1007/978-3-642-22061-6}}.

\end{thebibliography}
\bibliographystyle{elsarticle-num}

\end{document}